\documentclass[red,11pt,a4paper]{article}

\usepackage{cite}
\usepackage{amsmath}
\usepackage{amscd}
\usepackage{amssymb}
\usepackage{graphicx}
\usepackage{latexsym}
\usepackage{color}
\usepackage{verbatim}

\newcommand{\lr}[1]{\left( #1 \right)}
\newcommand{\lra}[1]{\langle #1 \rangle}
\newcommand{\lraa}[1]{\langle\langle  #1 \rangle \rangle}
\newcommand{\eq}[1]{\begin{equation}
\begin{split}
#1
\end{split}
\end{equation}}
\newcommand{\eqh}[1]{\begin{equation*}
\begin{split}
#1
\end{split}
\end{equation*}}
\newcommand{\Grad}{\nabla}
\newcommand{\Dt}{\frac{ d}{dt}}
\newcommand{\pt}{\partial_{t}}
\newcommand{\lap}{\Delta}
\newcommand{\pf}{{\noindent\it Proof.~}}
\newcommand{\ep}{\varepsilon}
\newcommand{\R}{\mathbb{R}}
\newtheorem{thm}{Theorem}
\newtheorem{lemma}[thm]{Lemma}
\newtheorem{prop}[thm]{Proposition}

\title{Kinetic theory of particle interactions mediated by dynamical networks}
\author{Julien Barr\'e$^{1,3}$, Pierre Degond$^2$, Ewelina Zatorska$^2$}

\begin{document}
\maketitle

\centerline{1.  Univ.  Nice Sophia Antipolis,  CNRS, Labo.  J.-A. Dieudonn\'e, UMR 7351,}
 \centerline{Parc Valrose, F-06108 Nice, France.}
 \bigskip
\centerline{2. Department of Mathematics, Imperial College London, }
\centerline{London SW7 2AZ, United Kingdom.}
\bigskip
\centerline{3. Institut Universitaire de France, 75005 Paris, France.}

\bigskip

\bigskip

\noindent{\bf Abstract:} We provide a detailed multiscale analysis of a system of particles interacting through a dynamical network of links. Starting from a microscopic model, via the mean field limit, we formally  derive coupled kinetic equations for the particle and link densities, following the approach of [{\it Degond et al., M3AS, 2016}].
Assuming that the process of remodelling the network is very fast, we simplify the description to a macroscopic model taking the form of single aggregation-diffusion equation for the density of particles. We analyze qualitatively this equation, addressing the stability of a homogeneous distribution of particles for a general potential.  For the Hookean potential we obtain a precise condition for the phase transition, and, using the central manifold reduction, we characterize the type of bifurcation  at the instability onset.

\bigskip


\section{Introduction}


 Cellular materials \cite{Gardel2004}, mucins \cite{BuCl07}, polymers \cite{BrDeYa, BaSu16} or social networks \cite{DeLiRi} are only few of the numerous examples of systems involving highly dynamical networks. A detailed modelling of these systems would require understanding complex chemical, biological or social  phenomena that are difficult to probe. Nevertheless, one common feature of these systems is the strong coupling between the  dynamical evolution of the individual agents (cells or monomers for instance) with that of the network mediating their interactions. The mathematical modelling of this strongly coupled dynamics is a challenging task, see for example \cite{PeDeDe}
but it is a necessary step towards building more complete models of complex biological or social phenomena.

The purpose of this paper is to provide a detailed multiscale analysis -- from a microscopic model to a macroscopic description, and its qualitative analysis -- of a system of particles interacting through a dynamical network, in a particularly simple setting: the basic entities are just point particles with local cross-links modelled by springs that are randomly created and destructed. 
In the mean field limit, assuming large number of particles and links as well as propagation of chaos, we  derive coupled kinetic equations for the particle and link densities. The link density distribution provides a statistical description of the network connectivity which turns out to be quite flexible and easily generalizable to other types of complex networks. See e.g. another application of this methodology to networks of interacting fibers in \cite{DeDePe}.

We focus on the regime where the network evolution triggered by the linking/unlinking processes happens on a very short timescale. 
In other words we are interested in observing dynamical networks on long time scale compared with the typical remodelling time scale.  In this regime the link density distribution becomes a local function of the particle distribution density. The latter evolves on the slow time scale through an effective equation which takes the form of an aggregation-diffusion equation, known also as the McKean-Vlasov equation \cite{McK67, ChPa}.
The applications of such an equation with different types of diffusion ranges from models of collective behavior of animals through granular media and chemotaxis models to self-assembly of nanoparticles, see \cite{SiSlTo15, KolCarBer, MoEd99, CaCaSc} and the references therein.
 In contrast to many of the aggregation-diffusion equations studied in the literature \cite{BeCaLa09, BeCh06, CaMcVi03, BeTo11} the model derived here features a compactly-supported potential. This model yields a very rich behavior, depending on two main parameters  describing the interaction range and the stiffness of the connecting links, that we investigate using both linear and nonlinear techniques. In particular, we identify the parameter ranges for the linear  stability/instability of the spatially homogeneous steady states. Moreover, the nonlinear analysis based on the central manifold reduction \cite{HaIo11} provides us with  a characterization of the type of bifurcation that appears at the  instability onset.  
Such bifurcations were previously studied in  \cite{ChPa} from a "thermodynamical" point of view, i.e. by looking at the minimizers of the free energy functional; we present here a dynamical point of view and make the connection with the thermodynamical approach.
In the case without diffusion, this free energy functional reduces to the interaction energy, whose minimizers have been studied in \cite{CaCaPa, SiSlTo15, CaDeMe};
for numerical studies in this direction we refer to  \cite{CaChHu}. In particular, global minimizers exist provided the associated potential is H-unstable, a classical notion in statistical mechanics linked to the phase transitions in the system \cite{FiRu66, Ru69}. Moreover, it was shown in \cite{CaCaPa}, that the minimizers are compactly supported for potentials with certain growth conditions at infinity. Generalization of these results to the case of compactly supported attraction-repulsion potential and linear diffusion, as in the system derived here, is a purpose of the future work.

The outline of the paper is the following. In the preliminaries of Section \ref{sec:model} we introduce an Individual-Based Model for  the point particles and the network, with rules for particles dynamics and network evolution. Then, in Section \ref{sec:kinetic}, we derive kinetic equations in a formal way following the approach from \cite{DeDePe} developed for systems of interacting fibers, when the number of particles $N$ and the number of links $K$ tend to infinity. In particular, we will assume that the ratio $K/N$ converges to some fixed positive limit $\xi$ that might be interpreted as an averaged number of links per particle. At the level of derivation of these equations, the precise character of particle interactions is not used and so the limit equations hold for a wide range of symmetric and integrable potentials.
In Section \ref{sec:macro}, we further simplify the description by assuming that the process of creating/destroying links is very fast. This enables us to derive a macroscopic model involving only the particle density, which takes the form of an aggregation-diffusion equation. In Section \ref{sec:genV}, we analyze qualitatively this macroscopic equation, addressing the stability of a homogeneous distribution of particles for a general potential, and in Section \ref{sec:HookeV} we address the same question for the Hookean potential, for which we obtain a precise condition for the bifurcation. Finally, in Section \ref{sec:non} we investigate via non linear analysis the character of the bifurcation, both for a rectangular (non degenerate unstable eigenvalue) and a square domain (degenerate unstable eigenvalue). In the last part of the paper, we illustrate
the criterion distinguishing between supercritical and subcritical bifurcations for the Hookean potential, and make connections with the very different approach by L. Chayes and V. Panferov in \cite{ChPa}.


\section{Modelling framework} \label{sec:model}


\subsection{Preliminaries}

The link between two particles located at the points $X_i$ and $X_j$ can be formed if their distance is less than a given radius of interaction $R$. If this condition is met the link is created in a Poisson process with probability $\nu_f^N$; it can be also destroyed with the probability $\nu_d^N$; both of them depend on $N$ -- the number of the particles in the whole system. When cross-linked, the particles interact with each-others subject to a pairwise potential 
\eq{\label{V}
V(X_i,X_j)= U(|X_i-X_j|).
}
For the moment we do not specify the character of interactions between the particles, trying to keep our derivation on a maximally general level.
 
We will first characterize the system of fixed number of particles, denoted by $N$, and fixed number of links, denoted by $K$. 
The equation of motion for each individual particle in the so-called overdamped regime, between two linking/unlinking events is:
\eq{\label{dX}
d X_i=-\mu \Grad_{X_i}Wdt+\sqrt{2D}dB_i,\quad i=1,\ldots,N.
}
Above, $B_i$ is a 2-dimensional Brownian motion $B_i=(B_i^1, B_i^2)$ with a positive diffusion coefficient $D>0$, $\mu>0$ is the mobility coefficient and  $W$ denotes the energy related to the maintenance of the links related to the potential $V$ as follows 
\eqh{
W=\sum_{k=1}^K V(X_{i(k)}, X_{j(k)}),}
where  $i(k), j(k)$ denote the indexes of particles connected by the link $k$. Plugging this definition into expression \eqref{dX}, we obtain
\eq{\label{dXi}
d X_i&=-{\mu}\sum_{k=1:i(k)=i}^K\left[\Grad_{x_1}V(X_{i(k)}, X_{j(k)})+\Grad_{x_2}V(X_{i(k)}, X_{j(k)})\right]dt+\sqrt{2D}dB_i\\
&= -{\mu}\sum_{k=1}^K\left[\delta_{i(k)}(i)\Grad_{x_1}V(X_{i(k)}, X_{j(k)})+\delta_{j(k)}(i)\Grad_{x_2}V(X_{i(k)}, X_{j(k)})\right]dt+\sqrt{2D}dB_i.
}
Our ultimate aim is to describe the systems of large number of particles. From the point of view of numerical simulations, the system of $N$ SDEs \eqref{dX} for large $N$, although fundamental, is too complex and thus costly to handle; it is also difficult to get a qualitative understanding of the behaviour of particles from \eqref{dX}. 
Therefore, in the next section we look for a "kinetic" description using probability distribution of particles and links rather then certain positions of each of the particles and links at a given time. 


\subsection{Derivation of the kinetic model}

\label{sec:kinetic}
We introduce the empirical distributions of the particles $f^N(x,t)$ and of the links $g^K(x_1,x_2,t)$, when the numbers of particles and links are finite and equal $N$ and $K$, respectively. They are equal to
\eqh{
&f^N(x,t)=\frac{1}{N}\sum_{i=1}^N \delta_{X_i}(x);\\
&g^K(x_1,x_2,t)=\frac{1}{2K}\sum_{k=1}^K\left[\delta_{X_{i(k)},X_{j(k)}}(x_1,x_2)+\delta_{X_{j(k)},X_{i(k)}}(x_1,x_2)\right],
}
where the symbol $\delta_{X_i}(x)$ is the Dirac delta centred at ${X_i(t)}$, with the similar definition for the two-point distribution. The above measures contain the full information about the positions of particles and links at time $t$.
For the sake of completeness we also introduce the two-particle empirical distribution 
\eq{h^N(x_1,x_2,t)=\frac{1}{N(N-1)}\sum_{i\neq j}\delta_{X_{i}(t),X_{j}(t)}(x_1,x_2).\label{def_h}}
Obviously, the two  distributions $h^N$ and $g^K$ are different, because not every pair of points is connected by a link.

The first part of this article is concerned with the derivation of the kinetic model obtained from \eqref{dX} in the mean-field limit. This process is roughly speaking a derivation of equations for the limit distributions $f$ and $g$, obtained from $f^N$ and $g^K$, by letting $N$ and $K$ to infinity, i.e.
\eqh{f(x,t):=\lim_{N\to\infty} f^N(x,t),\quad g(x_1,x_2,t)=\lim_{K\to \infty} g^K(x_1,x_2,t).
} 

The purpose of this section is to derive the equations for evolutions of particle and links distributions $f$ and $g$ in the limit of large number of particles and fibers. We have the following formal theorem. 
\begin{thm}\label{T1}
The kinetic system
\eq{\label{bdim}
&\pt f(x,t)={D\lap_x f(x,t)}+2\mu\xi\Grad_x\cdot F(x,t),\\
&\pt g(x_1,x_2,t)= {D\lr{\lap_{x_1} g(x_1,x_2,t)+\lap_{x_2}g(x_1,x_2,t)}}\\
&\hspace{2.5cm}+2\mu\xi \lr{  \Grad_{x_1}\cdot\lr{\frac{g(x_1,x_2)}{f(x_1)} F(x_1,t)}+\Grad_{x_2}\cdot\lr{\frac{g(x_1,x_2)}{f(x_2)} F(x_2,t)}}\\
&\hspace{2.5cm}+\frac{\nu_f}{2\xi} h(x_1,x_2,t)\chi_{|x_1-x_2|\leq R}-\nu_d g(x_1,x_2,t),
}
where
$$F(x,t)=\int g(x,y,t)\Grad_{x_1} V(x,y) dy,$$
and
\eqh{f(x,t):=\lim_{N\to\infty} f^N(x,t),\quad g(x_1,x_2,t)=\lim_{K\to \infty} g^K(x_1,x_2,t),\quad h(x_1,x_2,t)=\lim_{K\to \infty} h^N(x_1,x_2,t),
}
\eqh{
\nu_f= \lim_{N\to\infty} \nu_f^N(N-1),\quad \nu_d= \lim_{N\to\infty}\nu_d^N,
}
is a formal limit of the particle system \eqref{dX} as $N, K\to \infty$, provided that
\eqh{
\lim_{K,N\to\infty}\frac{K}{N}=\xi>0.}
\end{thm}
\pf The strategy of the proof is to first derive the equations for distribution of the particles $f^N(x,t)$ and of the links $g^K(x_1,x_2,t)$ in the situation when the number of each is finite and equal to $N$ and $K$. This happens between two linking/unlinking events in the time interval $(t,t+\Delta t)$. We will consider the behaviour of the system in this interval first and come back to the issue of creation of the new and destruction  of the old links in the end of the proof.

\bigskip

\noindent {\it Step 1.} Let us first introduce the notation that will allow us to identify both $f$ and $g$ with certain distributions. Following \cite{DeDePe} (Appendix A) we first introduce the one particle and two-particle observable functions, $\Phi(x)$ and $\Psi(x_1,x_2)$, respectively, and we define 
\eq{\label{fgf}
&\lra{ f^N(x,t),\Phi(x) }=\int f_1^N(x,t) \Phi(x) dx=\frac{1}{N}\sum_{i=1}^N \int \delta_{X_i(t)}(x) \Phi(x)\, dx=
\frac{1}{N}\sum_{i=1}^N \Phi(X_i(t));\\
&\lraa{g^K(x_1,x_2,t), \Psi(x_1,x_2)}=\int\int g^K(x_1,x_2,t) \Psi(x_1,x_2) dx_1\, dx_2\\
&\hspace{1,5cm}=\frac{1}{2K}\sum_{k=1}^K \int\int\left[\delta_{X_{i(k)},X_{j(k)}}(x_1,x_2)+\delta_{X_{j(k)},X_{i(k)}}(x_1,x_2)\right]\Psi(x_1,x_2) dx_1\, dx_2\\
&\hspace{1,5cm}=\frac{1}{2K}\sum_{k=1}^K \left[\Psi(X_{i(k)},X_{j(k)})+  \Psi(X_{j(k)},X_{i(k)})\right].
}
We will now apply the time derivative to the l.h.s. of these expressions and derive the equations of evolution of particles and links.

\bigskip

\noindent {\it Step 2.} We first derive the equation for the distribution of particles. Taking the time derivative of $\lra{ f^N(x,t),\Phi(x) }$  in \eqref{fgf} we get
\eqh{
\Dt\lra{ f^N(x,t),\Phi(x) }=
\frac{1}{N}\sum_{i=1}^N \Dt \Phi(X_i(t)). 
}
Using \eqref{dXi} and It\^{o}'s formula, we therefore obtain (formally)
\eq{\label{dist_f}
\Dt\lra{ f^N(x,t),\Phi(x) }&= -\mu \frac{1}{N}\sum_{i=1}^N \Grad_x\Phi(X_i(t))\cdot \nabla_{X_i}W\\
& \quad+D \frac{1}{N}\sum_{i=1}^N \Delta \Phi(X_i) +
\sqrt{2D} \frac{1}{N}\sum_{i=1}^N \Grad_x\Phi(X_i(t))\cdot \frac{dB_i}{dt}.
}
The random variables $\Grad_x\Phi(X_i(t))$ are not pairwise independent, since the $X_i$ are not independent. Nevertheless, the $dB_j$'s are pairwise independent, and are independent of $\Grad_x\Phi(X_i(t))$. Thus, the last term in \eqref{dist_f} is $1/N$ times the sum of uncorrelated random variables with zero expectation. Assuming, for instance, that the test functions have bounded derivatives, ensures that this last term is small in the $N\to \infty$ limit, so that it can be neglected in what follows. Thus
\eqh{
&\Dt\lra{ f^N(x,t),\Phi(x) }\\
&=-\frac{\mu}{N}\sum_{i=1}^N \Grad_x\Phi(X_i(t))\cdot 
\sum_{k=1}^K\left[\delta_{i(k)}(i)\Grad_{x_1}V(X_{i(k)}, X_{j(k)})+\delta_{j(k)}(i)\Grad_{x_2}V(X_{i(k)}, X_{j(k)})\right] \\
&\quad+D \frac{1}{N}\sum_{i=1}^N \Delta \Phi(X_i).
}
Exchanging the order of the sums with respect to $i$ and $k$ we get
\eq{\label{f1N}
&\Dt\lra{ f^N(x,t),\Phi(x) } \\
&\quad=-\frac{\mu}{N}\sum_{k=1}^K\left[\Grad_x\Phi(X_{i(k)})\cdot \Grad_{x_1}V(X_{i(k)}, X_{j(k)})+\Grad_x\Phi(X_{j(k)})\cdot \Grad_{x_2}V(X_{i(k)}, X_{j(k)})\right]\\
&\qquad+D \frac{1}{N}\sum_{i=1}^N \Delta \Phi(X_i)\\
&\quad=-\frac{\mu}{N}
\sum_{k=1}^K\left[\Grad_x\Phi(X_{i(k)})\cdot \Grad_{x_1}V(X_{i(k)}, X_{j(k)})+\Grad_x\Phi(X_{j(k)})\cdot \Grad_{x_1}V(X_{j(k)},X_{i(k)})\right]\\
&\qquad+D \frac{1}{N}\sum_{i=1}^N \Delta \Phi(X_i)\\
&\quad=-\frac{2K}{N}\lraa{g^K,\mu\Grad \Phi(x_1)\cdot\Grad_{x_1} V(x_1,x_2)} +D \lraa{f^N,\Delta \Phi}\\
&\quad=\frac{2\mu K}{N}\lraa{\Grad_{x_1}\cdot(g^K(x_1,x_2) \Grad_{x_1} V(x_1,x_2)), \Phi(x_1)} +D \lraa{\Delta f^N, \Phi}, 
}
where the second equality follows from the symmetry of potential $V$, i.e.  
$$\Grad_{x_2}V(X_{i(k)},X_{j(k)})=\Grad_{x_1}V(X_{j(k)},X_{i(k)});$$ 
the third equality in  \eqref{f1N} follows from the definition of distributions $f^N$ and $g^K$ \eqref{fgf}, and the last equality follows from integration by parts. Next, again formally, we can exchange the order of integration in \eqref{f1N}, so that it can be rewritten as
\eqh{\Dt
&\int f^N(x,t) \Phi(x) dx\\
&\quad= \frac{2\mu K}{N}\int\int \Grad_{x_1}\cdot(g^K(x_1,x_2) \Grad_{x_1} V(x_1,x_2)) \Phi(x_1) dx_1\, dx_2 +D\int \Delta f^N(x_1)\Phi(x_1)dx_1\\
&\quad=
\frac{2\mu K}{N}\int \Grad_{x_1}\cdot\lr{\int g^K(x_1,x_2) \Grad_{x_1} V(x_1,x_2)dx_2} \Phi(x_1) dx_1 +D\int \Delta f^N(x_1)\Phi(x_1)dx_1.
}
Therefore, letting $N,K$ to infinity, assuming that $\frac{K}{N}\to\xi$ and that there exist the limits 
\eqh{\lim_{N\to\infty}f^N=f\quad\mbox{and}\quad \lim_{K\to\infty}g^K=g,}
 we obtain (after change of variables $x_1\to x$, $x_2\to x'$) a distributional formulation of equation for $f$. The differential form of this equation is
\eq{\label{ev_f1}
{\pt f(x,t)=2\mu\xi\Grad_x\cdot F(x,t) +D \Delta f,\qquad
F(x_1,t)=\int g(x_1,x_2,t)\Grad_{x_1} V(x_1,x_2) dx_2.}
}

\bigskip

\noindent {\it Step 3.}
After deriving the equation for distribution of particles $f$ we want to derive the equation for $g$ in the analogous way.
We remark that the noise in \eqref{dXi} transforms directly into a linear diffusion term for $f$, all other contributions  vanish in the large $N$ limit. It is not difficult to see that the same simplification takes place for $g^K$ in the $K\to \infty$ limit. Thus, to reduce the computations we will first use \eqref{dXi} without noise, and reintroduce the diffusion term in the end.

Taking the time derivative of the second equality in \eqref{fgf} we obtain
\eq{\label{gE1E2}
&\Dt\lraa{g^K(x_1,x_2,t), \Psi(x_1,x_2)} \\
&=\frac{1}{2K}\sum_{k=1}^K \left[\Grad_{x_1}\Psi(X_{i(k)},X_{j(k)})\cdot \Dt X_{i(k)}+  \Grad_{x_1}\Psi(X_{j(k)},X_{i(k)})\cdot \Dt X_{j(k)}\right]\\
&\hspace{0,2cm}+\frac{1}{2K}\sum_{k=1}^K\left[\Grad_{x_2}\Psi(X_{i(k)},X_{j(k)})\cdot \Dt X_{j(k)}+  \Grad_{x_2}\Psi(X_{j(k)},X_{i(k)})\cdot \Dt X_{i(k)}\right]\\
&=E_1+E_2.
}
We now present how to treat $E_1$, $E_2$ can be handled analogously. We first use \eqref{dXi} (without noise) to write
\eq{\label{E1}
E_1=&\frac{-\mu}{2K}\sum_{k=1}^K\bigg\{\Grad_{x_1}\Psi(X_{i(k)},X_{j(k)})\\
&\hspace{1,3cm}\times \sum_{k'=1}^K\left[\delta_{i(k')}(i(k))\Grad_{x_1}V(X_{i(k')}, X_{j(k')})+\delta_{j(k')}(i(k))\Grad_{x_2}V(X_{i(k')}, X_{j(k')})\right]\bigg\}
 \\
&\frac{-\mu}{2K}\sum_{k=1}^K\bigg\{\Grad_{x_1}\Psi(X_{j(k)},X_{i(k)})\\
&\hspace{1,3cm}\times \sum_{k'=1}^K\left[\delta_{i(k')}(j(k))\Grad_{x_1}V(X_{i(k')}, X_{j(k')})+\delta_{j(k')}(j(k))\Grad_{x_2}V(X_{i(k')}, X_{j(k')})\right]\bigg\}
 \\
 =&\frac{-\mu}{2K}\sum_{k=1}^K\bigg\{\Grad_{x_1}\Psi(X_{i(k)},X_{j(k)})\\
&\hspace{1,3cm}\times\sum_{k'=1}^K\left[\delta_{i(k')}(i(k))\Grad_{x_1}V(X_{i(k')}, X_{j(k')})+\delta_{j(k')}(i(k))\Grad_{x_1}V(X_{j(k')}, X_{i(k')})\right]\bigg\}
 \\
&\frac{-\mu}{2K}\sum_{k=1}^K\bigg\{\Grad_{x_1}\Psi(X_{j(k)},X_{i(k)})\\
&\hspace{1,3cm}\times \sum_{k'=1}^K\left[\delta_{i(k')}(j(k))\Grad_{x_1}V(X_{i(k')}, X_{j(k')})+\delta_{j(k')}(j(k))\Grad_{x_1}V(X_{j(k')}, X_{i(k')})\right]\bigg\}
 \\
=&\frac{-\mu}{2K}\sum_{k'=1}^K \bigg\{\Grad_{x_1}V(X_{i(k')}, X_{j(k')})\\
&\hspace{1,3cm}\times  
\sum_{k=1}^K\left[ \delta_{i(k')}(i(k))\Grad_{x_1}\Psi(X_{i(k)},X_{j(k)})+\delta_{i(k')}(j(k))\Grad_{x_1}\Psi(X_{j(k)},X_{i(k)})\right]
 \bigg\}\\
&\frac{-\mu}{2K}\sum_{k'=1}^K \bigg\{\Grad_{x_1}V(X_{j(k')}, X_{i(k')})\\
&\hspace{1,3cm}\times 
\sum_{k=1}^K\left[ \delta_{j(k')}(i(k))\Grad_{x_1}\Psi(X_{i(k)},X_{j(k)})+\delta_{j(k')}(j(k))\Grad_{x_1}\Psi(X_{j(k)},X_{i(k)})\right]\bigg\}
}
We see that the first sum with respect to $k$ in the last equality of \eqref{E1}, i.e.
\eq{\label{sumE1}
\sum_{k=1}^K\left[ \delta_{i(k')}(i(k))\Grad_{x_1}\Psi(X_{i(k)},X_{j(k)})+\delta_{i(k')}(j(k))\Grad_{x_1}\Psi(X_{j(k)},X_{i(k)})\right]
}
does not vanish if
either $i(k)=i(k')$ or $j(k)=i(k')$. To understand it better let us look at the link number $k'$. Its beginning is $i(k')$ and it is a certain fixed particle as was the link.

If we now compute the above sum neglecting the Kronecker symbols we get 2K of different elements. But for the Kronecker symbols included we act in the following way: we take the first link $k=1$ and check if 
$i(1)=i(k')$ if yes then definitely $j(1)\neq i(k')$ thus the first element of the sum is equal to $\Grad_{x_1}\Psi(X_{i(1)},X_{j(1)})$, if $i(1)\neq i(k')$ then we check if $j(1)=i(k')$ if yes the first element of the sum equals 
$\Grad_{x_1}\Psi(X_{j(1)},X_{i(1)})$. Finally if $i(k')\neq i(1)$ and $i(k')\neq j(1)$ the above sum reduces to the subset $k\geq 2$.
Hence the maximal number of elements of the above sum is $K$, but in fact it will be equal to the number of links connected to $i(k')$ and it may be less then the number of all links  $K$.

We now introduce a number of links connected to $i(k')$
\eqh{
C_{i(k')}=\#\{k\ |\ i(k)=i(k')\ \ \mbox{or}\ \ j(k)=i(k')\}.
}
Thus, dividing \eqref{sumE1} by $C_{i(k')}$ and letting $K\to\infty$ gives rise to a certain probability associated with $i(k')$, we have
\eq{\label{Ci}
\lim_{K\to \infty}\frac{1}{C_{i(k')}}\sum_{k=1}^K\left[ \delta_{i(k')}(i(k))\Grad_{x_1}\Psi(X_{i(k)},X_{j(k)})+\delta_{i(k')}(j(k))\Grad_{x_1}\Psi(X_{j(k)},X_{i(k)})\right]\\
=2\int{(\Grad_{x_1}\Psi P)(X_{i(k')},x_2)}dx_2,
}
where 
\eqh{P(X_{i(k')},x_2)=\frac{g(X_{i(k')},x_2)}{\int g(X_{i(k')},x_2) dx_2}}
is a conditional probability of finding a link, provided one of its ends is at $X_{i(k')}$.

We can now estimate the limit of mean number of links per particle when $N,K\to\infty$, $\frac{K}{N}\to\xi$. Around the point $X_{i(k')}$ we have
\eqh{C_{i(k')}=\frac{K\int g^K(X_{i(k')},x_2) dx_2}{N f^N(X_{i(k')})},}
therefore
\eq{\label{Cilim}
\lim_{K,N\to \infty,\, \frac{K}{N}\to\xi} C_{i(k')}= \xi \frac{\int g(X_{i(k')},x_2) dx_2}{f(X_{i(k')})}.}
Combining \eqref{Ci} and \eqref{Cilim}, we obtain
\eqh{
\lim_{N,K\to \infty,\, \frac{K}{N}\to\xi}&\sum_{k=1}^K\left[ \delta_{i(k'),j(k)}\Grad_{x_1}\Psi(X_{i(k)},X_{j(k)})+\delta_{i(k'),i(k)}\Grad_{x_1}\Psi(X_{j(k)},X_{i(k)})\right]\\
&=\frac{2\xi}{f(X_{i(k')})}\int{(\Grad_{x_1}\Psi g)(X_{i(k')},x_2)}dx_2,
}
thus the limit of \eqref{E1} reads
\eqh{
\lim_{K,N\to \infty,\, \frac{K}{N}\to\xi} E_1 &  
=\lim_{K\to\infty} -\frac{\mu\xi}{K}\sum_{k'=1}^K \left[\Grad_{x_1}V(X_{i(k')}, X_{j(k')})
  \cdot \frac{\int{(\Grad_{x_1}\Psi g)(X_{i(k')},x_2)}dx_2}{f(X_{i(k')})}
 \right.\\
&\hspace{2,5cm} \left.+
\Grad_{x_1}V(X_{j(k')}, X_{i(k')}) \cdot
\frac{\int{(\Grad_{x_1}\Psi g)(X_{j(k')},x_2)}dx_2}{f(X_{j(k')})}
 \right]\\
 &=-2\mu\xi\lraa{g,\Grad_{x_1}V(x_1,x_2) \cdot
\frac{\int{(\Grad_{x_1}\Psi g)(x_1,x_2)}dx_2}{f(x_1)}}
.}
Now, coming back to \eqref{gE1E2} and performing the same procedure for $E_2$ we obtain
\eqh{
\Dt\lraa{g(x_1,x_2,t), \Psi(x_1,x_2)} &=-2\mu\xi\lraa{g,\Grad_{x_1}V(x_1,x_2) \cdot
\frac{\int{(\Grad_{x_1}\Psi g)(x_1,x_2)}dx_2}{f(x_1)}}\\
&\quad-2\mu\xi\lraa{g,\Grad_{x_1}V(x_1,x_2) \cdot
\frac{\int{(\Grad_{x_2}\Psi g)(x_2,x_1)}dx_2}{f(x_1)}}.
}
Integrating by parts, changing the variables and order of integrals we easily obtain
\eqh{
\Dt\lraa{g(x_1,x_2,t), \Psi(x_1,x_2)} &=2\mu\xi\lraa{\Grad_{x_1}\cdot\lr{\frac{g(x_1,x_2)}{f(x_1)}\int g\Grad_{x_1}V(x_1,x_2) dx_2}, \Psi(x_1,x_2)}\\
&\quad+2\mu\xi\lraa{\Grad_{x_2}\cdot\lr{\frac{g(x_1,x_2)}{f(x_2)}\int g\Grad_{x_1}V(x_2,x_1) dx_1}, \Psi(x_1,x_2)}.
}
Therefore, the differential form of equation for $g$ reads
\eq{\label{ev_g}
\pt g(x_1,x_2,t)=& D\lr{\lap_{x_1} g(x_1,x_2,t)+\lap_{x_2}g(x_1,x_2,t)}\\
&+2\mu\xi \Grad_{x_1}\cdot\lr{\frac{g(x_1,x_2)}{f(x_1)} F(x_1,t)}+2\mu\xi \Grad_{x_2}\cdot\lr{\frac{g(x_1,x_2)}{f(x_2)} F(x_2,t)},
}
where we have reintroduced the diffusion terms due to the noise in \eqref{dXi}, and $F(x_1)$ is the same one as defined as in \eqref{ev_f1}, recall
\eqh{
F(x_1,t)=\int g(x_1,x_2,t)\Grad_{x_1} V(x_1,x_2) dx_2,\quad F(x_2,t)=\int g(x_2,x_1,t)\Grad_{x_1} V(x_2,x_1) dx_1.}

\noindent {\it Step 4.} Equations \eqref{ev_f1} and \eqref{ev_g} do not take into account the phenomena of creation and destruction of links. According to the description at the beginning of this paper, our model describes a process of creation of links with the probability $\nu_f^N$, provided the two particles are sufficiently close to each others. Surely, the number of new links  will be proportional to the number of couples of the particles such that one of them is close to $x_1$ and the other one is close to $x_2$, whose distance is less than $R$, this number is equal to:
$$\frac{N(N-1)}{2} h(x_1,x_2,t)\chi_{|x_1-x_2|\leq R}\, dx_1\,dx_2\, dt ,$$
where $h(x_1,x_2,t)=\lim_{N\to \infty} h^N$ and $h^N=h^N(x_1,x_2,t)$ is the two-particle distribution defined in \eqref{def_h}.
This number has to be decreased by the number of couples that are already connected by existing links:
$$K g(x_1,x_2,t)\,dx_1\,dx_2\, dt.$$
Therefore, the number of the new links created during the time interval $[t,t+dt[$
between two points $x_1$ and $x_2$ is equal to 
$$\nu_f^N \lr{\frac{N(N-1)}{2}h(x_1,x_2,t)\chi_{|x_1-x_2|\leq R}-K g(x_1,x_2,t)
}\, dx_1\,dx_2\, dt .$$
Dividing this expression by $K$ used for normalization of function $g$ 
and letting $N,K\to \infty$ so that $\frac{K}{N}\to\xi$ and $\nu_f^N(N-1)\to  \nu_f$ we obtain the probability of creation of the new link equal to
$$\frac{ \nu_f}{2\xi}h(x_1,x_2,t)\chi_{|x_1-x_2|\leq R}.$$
Similarly, the probability that the existing link will be destroyed in the same time interval $[t,t+dt[$ is equal to
$$\nu_d g(x_1,x_2,t),$$
where we used $\nu_d=\lim_{N\to\infty}\nu_d^N$.
If we now include these source terms in \eqref{ev_g}, we get
\eqh{
&\pt g(x_1,x_2,t)= {D\lr{\lap_{x_1} g(x_1,x_2,t)+\lap_{x_2}g(x_1,x_2,t)}}\\
&\hspace{2.5cm}+2\mu\xi \lr{  \Grad_{x_1}\cdot\lr{\frac{g(x_1,x_2)}{f(x_1)} F(x_1,t)}+\Grad_{x_2}\cdot\lr{\frac{g(x_1,x_2)}{f(x_2)} F(x_2,t)}}\\
&\hspace{2.5cm}+\frac{\nu_f}{2\xi} h(x_1,x_2,t)\chi_{|x_1-x_2|\leq R}-\nu_d g(x_1,x_2,t).
}
This, together with equation \eqref{ev_f1} gives the system \eqref{bdim}. Theorem \ref{T1} is proved. $\Box$

Note that system \eqref{bdim} is not closed, since all the three distributions $f,\ g$ and $h$ are a-priori unknown. In order to close this system we will have to introduce some closure assumption; this will be done in the next section.


\section{Derivation of the macroscopic equations}\label{sec:macro}


The equations of distributions of particles and links in the form introduced in Theorem \ref{T1} do not reveal anything more than relations between certain mechanisms leading to evolution in time of $f$ and $g$. To get somehow deeper insight to the behaviour of the system we introduce the characteristic values of the physical quantities appearing in the system. We denote by $t_{0}$ the unit of time and by $x_{0}$ the unit of space. Accordingly, we also identify the units for the parameters of the system and their dimensionless values
\eqh{\mu'=\frac{\mu}{t_{0}},\quad D'=\frac{Dt_{0}}{x_{0}^2}, \quad\nu_f'=\nu_f t_0,\quad \nu_d'=\nu_d t_0}
The units of distribution functions are
$$f_{0}=\frac{1}{x_{0}^2},\qquad g_{0}=\frac{1}{x_{0}^4}, \qquad h_{0}=\frac{1}{x_{0}^4},$$
where the powers reflect the fact that the physical domain is two-dimensional and the dimensionless values are given by $X'=X/X_{0}$, therefore
\eqh{
\pt f(t,x)=\frac{1}{t_0x_0^2}\partial_{t'}f'(t',x'), \quad \pt g(t,x)=\frac{1}{t_0x_0^4}\partial_{t'}g'(t',x'),\quad \pt h(t,x)=\frac{1}{t_0x_0^4}\partial_{t'}h'(t',x').
}
Similarly if we assume that the potential scales as the potential energy $V_{0}=\frac{x_{0}^2}{t_0^2}$, thus
\eqh{
\Grad_{x_1} V(x_1,x_2)&=  \frac{x_{0}}{t_{0}^2} \Grad_{x_1'} V'(x_1', x_2'),\\
 \Grad_{x_1}\cdot F(x_1)&=\Grad_{x_1}\cdot\int g(x_1,x_2)\Grad_{x_1} V(x_1,x_2) dx_2\\
 &=
x_0^{-1}\Grad_{x_1'}\cdot\int x_0^{-4}\tilde g(\tilde x_1,\tilde x_2) \frac{x_{0}}{t_{0}^2} \Grad_{x_1'} V'(x_1', x_2') x_0^2d x_2'\\
&= \frac{1}{x_0^2t_0^2} \Grad_{x_1'}\cdot F'( x_1'),
}
Substituting the above formulas into \eqref{ev_f1} and omitting the primes, we obtain
the scaled version of equation for $f$:
\eq{\label{ev_f1p}
{\pt f(x,t)=2\mu\xi \Grad_x\cdot F(x,t) +D  \Delta f,}
}
and the scaled version of equation for $g$:
\eqh{
&\pt g(x_1,x_2,t)= {D\lr{\lap_{x_1} g(x_1,x_2,t)+\lap_{x_2}g(x_1,x_2,t)}}\\
&\hspace{2.5cm}+2\mu\xi  \lr{  \Grad_{x_1}\cdot\lr{\frac{g(x_1,x_2)}{f(x_1)} F(x_1,t)}+\Grad_{x_2}\cdot\lr{\frac{g(x_1,x_2)}{f(x_2)} F(x_2,t)}}\\
&\hspace{2.5cm}+\frac{\nu_f}{2\xi} h(x_1,x_2,t)\chi_{|x_1-x_2|\leq R}-\nu_d g(x_1,x_2,t).
}
These equations give us some freedom in the choice of the time scale and space scale, from now on we will use time and space units such that 
$$\mu=1\quad \mbox{and} \quad D=1.$$
The next step is to introduce the macroscopic scaling for these units using small parameter $\ep<<1$: $ x''_0=\ep^{-1/2}x_0$, $ t''_0=\ep^{-1}t_0$.
Then the new variables and unknowns are
\eqh{
x''&=\ep^{1/2}x,\quad  t''=\ep t, \quad f''(x'')=\ep^{-1} f(x),\\
g''(x''_1, x''_2)&=\ep^{-2}g(x_1,x_2),\quad h''( x''_1,x''_2)=\ep^{-2}h(x_1,x_2).}
Then, we also introduce the scaling of the potential \eqref{V}. This time, we assume a small intensity of interactions, therefore $V(x_1,x_2)\approx V''( x''_1, x''_2)$,
moreover,
\eqh{
\Grad_x V(x_1,x_2)&=\ep^{1/2} \Grad_{x''}  V''(x''_1,x''_2),\\
\Grad_{x_1} F(x_1)&=\Grad_{x_1}\int g(x_1,x_2)\Grad_{x_1} V(x_1,x_2) dx_2\\
&=
\ep^{1/2}\Grad_{\tilde x_1}\int \ep^2  g''( x''_1,x''_2)\ep^{1/2}\Grad_{x''_1} V''( x''_1,x''_2) \ep^{-1} dx''_2\\
&= \ep^{2} \Grad_{x''_1} F''( x''_1),
}
so when we compare the terms of order $\ep^{2}$ in expansion of $f$ in \eqref{ev_f1p} with $\mu, D=1$, we basically get the same equation for $f''$
\eq{\label{eq:f}
\partial_{ t''}  f''={\lap_{ x''} f''}+2\xi\Grad_{x''}\cdot F''.}
Concerning the equation for distribution of links, we assume that the creation and destruction of links is a very fast process, meaning that
\[\nu''_f=\ep^2\nu_f,\quad \nu''_d=\ep^2\nu_d,\]
 noticing that $\chi_{|x_1-x_2|\leq R}=\chi_{| x''_1- x''_2|\leq R''}$, we have
 \eq{\label{eq:g}
 &\ep^3 \partial_{ t''}  g''\\
 &\qquad =\ep^3\lap g''
+2\xi \left[ \ep^{1/2} \Grad_{x''_1}\cdot\lr{\frac{\ep^2 g''}{\ep f''(x''_1)}\ep^{3/2}  F''(x''_1)}+\ep^{1/2}\Grad_{x''_2}\cdot\lr{\frac{\ep^2 g''}{\ep f''(x''_2)}\ep^{3/2}  F''_1(x''_2)}\right]\\
&\hspace{2.5cm}+\ep^2\lr{\frac{\nu_f}{2\xi} h''\chi_{|x''_1- x''_2|\leq  R''}-\nu_d  g''}\\
&\qquad =\ep^3\lr{\lap g''
+2\xi \left[  \Grad_{ x''_1}\cdot\lr{\frac{g''}{ f''(x''_1)} F''( x''_1)}+\Grad_{x''_2}\cdot\lr{\frac{g''}{f''(x''_2)}  F'(x''_2)}\right]}\\
&\hspace{2.5cm}+\lr{\frac{\nu''_f}{2\xi}  h''\chi_{|x''_1- x''_2|\leq R''}-\nu''_d g''}.
 }
 Our purpose now is to let $\ep$ to zero in \eqref{eq:f} and \eqref{eq:g}. Assuming again that $ f''$, $ g''$ and $ h''$ exist we denote $f_\ep= f''$, $g_\ep=g''$, $h_\ep= h''$, we then have the following proposition.
 \begin{prop}
 Assume that $h_\ep(x_1,x_2)=f_\ep(x_1)f_\ep(x_2)$,  and that $V(X_i,X_j)= U(|X_i-X_j|)$, then provided the following limits exist
 $$f:=\lim_{\ep\to0} f_\ep,\quad g:=\lim_{\ep\to0} g_\ep$$
 they formally satisfy
 \begin{subequations}\label{main_sys}
\eq{
\partial_t f(t,x)=\lap_x f(t,x)+ 
\frac{\nu_f}{\nu_d}\Grad_x\cdot (f(t,x)\Grad_x (\tilde V\ast f)(t,x)) \label{main_sys1}}
\eq{g(t,x,y)=\frac{\nu_f}{2\xi \nu_d}  f(t,x) f(t,y)\chi_{|x-y|\leq R},\label{main_sys2}}
\end{subequations}
for some compactly supported potential $\tilde V$ specified below.
 \end{prop}
 \pf Let us start with the limit equation for the distribution of links. From \eqref{eq:g}, using the assumption on small correlations we obtain 
\eqh{\frac{\nu_f}{2\xi} f_\ep( t,x) f_\ep(t, y)\chi_{|x-y|\leq R}-\nu_d g_\ep(t,x,y)= O(\ep^3).
}
Letting $\ep\to0$ in the above formula, we formally obtain  \eqref{main_sys2}, which is an explicit formula for $g$.
Therefore, plugging this relation into \eqref{eq:f} and dropping the tildes again we obtain the equation for $f$:
\eqh{
\partial_{t}  f={\lap_{x}  f}+\Grad_{ x}\cdot  F,\qquad F=\frac{\nu_f}{\nu_d} f(x)
\int   f(y)\Grad_{x}V(x,y) \chi_{|x-y|\leq R}dy.
}
Taking into account the form of the potential, we can rewrite the above equation in slightly different form
\eq{\label{tV}
\partial_t f=\lap_x f+ \frac{\nu_f}{\nu_g}\Grad_x\cdot \lr{f(x)\int \Grad\tilde V(x-y) f(y) dy  }
}
for some $\tilde V$ such that
$$\Grad_i \tilde V(x)= U'(|x|)\chi_{|x|\leq R}\vec{e}_i,\quad i=1,2,$$
which gives \eqref{main_sys1}. $\Box$


\section{Analysis of the macroscopic equation: general potential}\label{sec:genV}


\subsection{Remark about the free energy}

The above system, particularly equation \eqref{main_sys1}, is well known in the literature as an aggregation-diffusion equation, also as McKean-Vlasov equation. For analytical and numerical results devoted to solvability and asymptotic analysis of solutions, depending on the shape of the potential $\tilde V$, see for instance \cite{ChPa, CaCaSc}.
Concerning the steady states,  an exhaustive analysis of this problem would require finding the minima of the following energy functional associated with \eqref{main_sys1}:
\eq{
{\cal F}(f)=\int \lr{ f\log f+\frac{1}{2} \frac{\nu_f}{\nu_d}f(\tilde V\ast f)} dx.
\label{eq:entropy}
}
It is easy to check that ${\cal F}(t)$ is dissipated in time:
\eqh{
&\Dt {\cal F}(f)=\int \lr{\pt f\log f+\pt f+\frac{\nu_f}{\nu_d}\pt f(\tilde V\ast f) }dx\\
&=\int \bigg(\lap f\log f+  \frac{\nu_f}{\nu_d}\Grad\cdot(f\Grad(\tilde V\ast f))\log f \\
&\hspace{2,3cm}+\frac{\nu_f}{\nu_d}(\tilde V\ast f)\lap f+\lr{\frac{\nu_f}{\nu_d}}^2\Grad\cdot(f\Grad (\tilde V\ast f))(\tilde V\ast f) \bigg)dx\\
&=-\int \lr{\frac{|\Grad f|^2}{f}+ \lr{2 \frac{\nu_f}{\nu_d}}\Grad (\tilde V\ast f)\cdot\Grad f +\lr{ \frac{\nu_f}{\nu_d}}^2f|\Grad (\tilde V\ast f)|^2 }dx\\
&=- \int\lr{\frac{\Grad f}{f^{1/2}}+  \frac{\nu_f}{\nu_d}f^{1/2}\Grad (\tilde V\ast f)}^2dx
\leq 0.
}


\subsection{Constant steady states}

In this note, we want to focus only on the constant steady states, i.e. $f_\star={\rm const}$, which, on bounded domains, have an interpretation as probability measures. It turns out that the stability or instability of the steady states for \eqref{main_sys1} is related to the notion of H-stability of the potential $\tilde V$. According to the definitions from classical statistical mechanics, the compactly supported potential $\tilde V$ is H-stable provided the integral  $\int_{\R^2}  \tilde V(x)\, dx$ is positive, otherwise it is not H-stable (unstable) \cite{Ru69}. For the H-stable potentials, the aggregation part of equation  \eqref{main_sys1} acts as diffusion, so, any initial perturbation is smoothen infinitely fast. For potentials that are not H-stable, the asymptotical behaviour of the solution is much more interesting. For our system in its general form we only prove the following criterion for instability of the constant steady states.
\begin{lemma}\label{L:un}
Let the potential $\tilde V$ be integrable and let
\eqh{M= \int_{\R^2}  \tilde V(x)\, dx<0.
}
Then the constant steady state $f_\star$ is unstable  if
\eq{f_\star>\frac{-1}{ M}\frac{\nu_d}{\nu_f}.\label{instab}}
\end{lemma}
\pf In order to check the stability of the constant steady state $f_\star>0$, we linearize \eqref{main_sys1} around $f_\star$. We assume that $f$ is a small perturbation of $f_\star$ ($f<< f_\star$) and thus $f$ satisfies
\eq{\label{main_per}
\partial_t f(t,x)=\lap_x f(t,x)+ 
f_\star\frac{\nu_f}{\nu_d}\lap_x( (\tilde V\ast f)(t,x)).
}
Then  we  apply the Fourier transform in space to both sides of \eqref{main_per}, we obtain
\eq{\label{tr_f}
\pt\hat f(t,y)=-y^2\hat f(t,y)-2\pi f_\star\frac{\nu_f}{\nu_d} y^2\hat{\tilde V}\hat f(t,y).
}
The Taylor expansion around zero of the Fourier transform of $\tilde V$ is equal to
\eqh{
\hat{ \tilde V}(y)=\frac{1}{2\pi}
\int_{\R^2}  \tilde V(x)\, dx+ O(y)= \frac{M}{2\pi} + O(y).
}
Plugging it into \eqref{tr_f} we obtain
\eq{\label{fur_U}
\pt\log\hat f(t,y)=-\lr{1+ f_\star\frac{\nu_f}{\nu_d} M} y^2+ O(y^3),
}
and so, for negative $M$, we can always find sufficiently large  $f_\star$ leading to instability of the steady state $f_\star$. More precisely, for \eqref{instab}
the r.h.s. of \eqref{fur_U} for sufficiently small $y$ is larger then some positive constant $c$, thus
$$\hat f(t)\geq\hat f_0 e^{ct}\to\infty\quad \mbox{for}\quad {t\to\infty},$$
and so, the steady state $f_\star$ is unstable. $\Box$


\section{Analysis of the macroscopic equation: Hookean potential}\label{sec:HookeV}


\subsection{Preliminaries}

Until this moment, the exact form of potential \eqref{V} did not play any role and we could work assuming only its symmetricity  and integrability. Let us now focus on a particular form. If we imagine that the links between the particles act like springs, the interaction potential is given by the Hooke law
\eqh{V(x_1,x_2)=\frac{\kappa}{2}\lr{|x_1-x_2|- l_0}^2,}
where $l_0$ denotes the rest length of the spring and the intensity parameter $\kappa$ is a positive number, characteristic of the spring. We then have
\eqh{\int  f(y)\chi_{|y- x|\leq R}\Grad_{x}V(x,y) dy=\int f(y) \kappa(|x-y|-l_0)\frac{x-y}{|x-y|}\chi_{|x-y|\leq R} dy. 
}
We now want to find $\tilde V$ such that the equation for $f$ is in the form \eqref{tV}.
In our case $\tilde V(x)$ satisfies $\Grad_{i}\tilde V(x)=\kappa(|x|-l_0)\chi_{|x|\leq R}\vec{e}_i,$ where $x\in \R^2$,  moreover $\tilde V(x)=0$ for $|x|>R$.
First, it is easy to see that $\tilde V(x)$ is a radially symmetric function, thus we can introduce $U(|x|)=\tilde V(x)$, secondly since the potential $U(r)$ vanishes for $r\geq R$ we have
\eqh{U(2R)-U(r)=\int_r^R (s-l_0) ds =\frac{\kappa}{2}\left[(R-l_0)^2-(r-l_0)^2\right].}
Therefore, $U(r)=\frac{\kappa}{2}\left[(r-l_0)^2-(R-l_0)^2\right]$, and so
\eq{\label{TVH}
\tilde V(x)=
\left\{
\begin{array}{lll}
\frac{\kappa}{2}\left[(|x|-l_0)^2-(R-l_0)^2\right],& \mbox{for}&|x|< R,\\
0& \mbox{for}&|x|\geq R,
\end{array}
\right.
}
see the picture below.
\begin{center}
\includegraphics[scale=0.35]{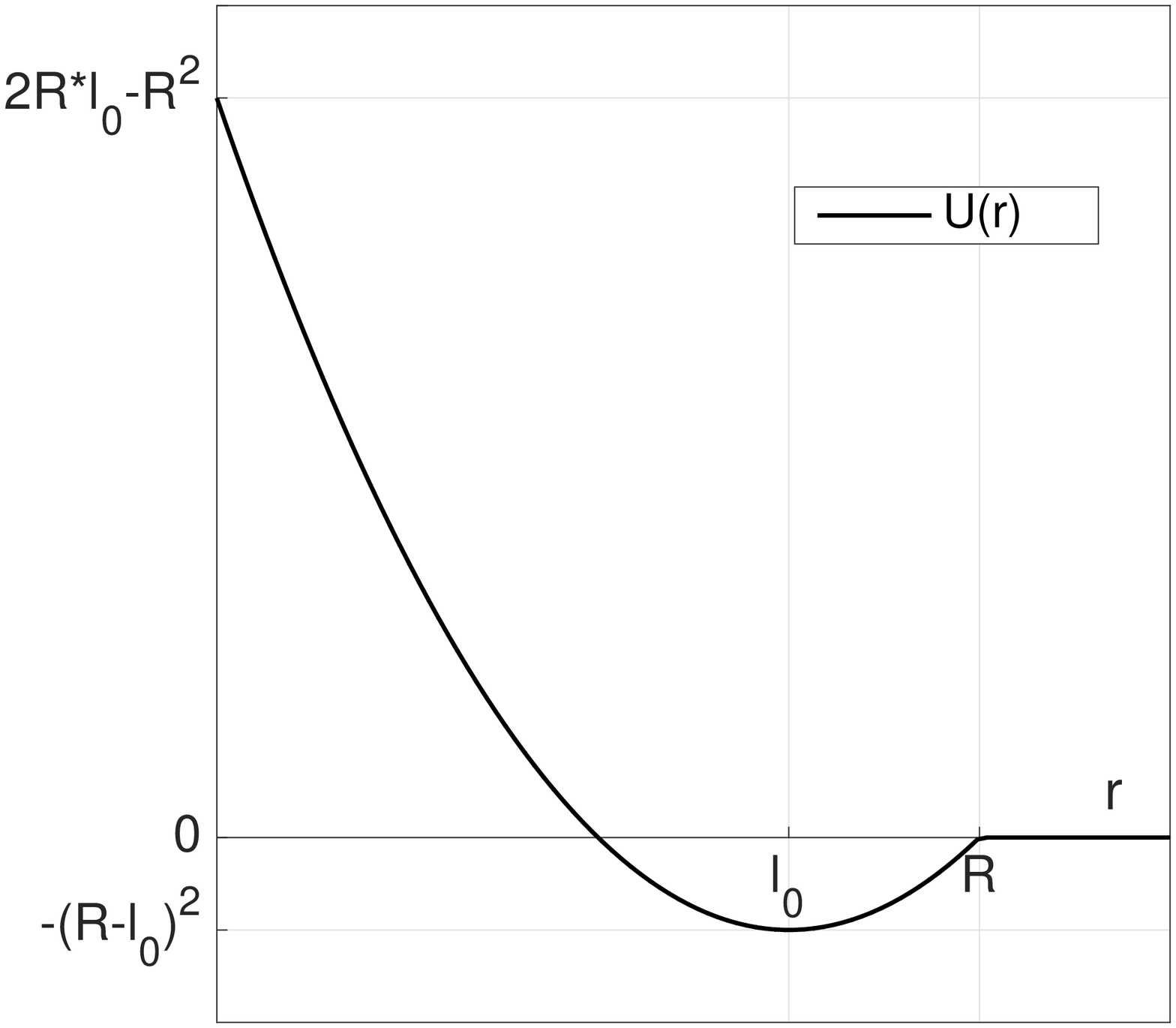}
\end{center}
Let us now compute the integral of our potential $\tilde V$ given in \eqref{TVH}. We have
 \eqh{
\int_{\R^2}{\tilde V(x) dx}&=\frac{\kappa}{2}\int_{\R^2}{\left[(|x|-l_0)^2-(R-l_0)^2\right] \chi_{|x|<R}\,dx}
=\pi\kappa\int_{0}^R{\left[(r-l_0)^2-(R-l_0)^2\right] r\, dr}\\
&=\pi\kappa\lr{\frac{r^4}{4}-\frac{2r^3 l_0}{3}-\frac{R^2r^2}{2}+R l_0 r^2}\Big|_{0}^{R}=\pi\kappa R^3\lr{\frac{ l_0}{3}-\frac{R}{4}},
}
therefore, according to the definition given above, $\tilde V$ is H-stable if the condition $l_0>\frac{3R}{4}$ is satisfied. Lemma \ref{L:un}
provided a special criterion for the constant steady state to be unstable, and this is basically all the information we can get for the whole space case. However, if we now consider the same problem on the space periodic domain the criteria obtained in Lemma \ref{L:un} will have to include the size of the domain. Moreover, it can happen that even if unstable, the steady state might be only weakly unstable, meaning that only one mode from countable set of modes will be unstable, while the rest of them will be stable.
The intention of the linear analysis in the whole space case presented below is to provide some intuition on the behaviour of the potential, so that it is more intuitive how to "select" the unstable modes in the second part of this section.


\subsection{Linear analysis in the whole space}

To understand the behaviour of the solutions close to the stability/instability threshold  \eqref{instab} we come back to equation \eqref{tr_f} and we compute the Fourier transform of $\tilde V$ given by \eqref{TVH}
\eqh{
\hat{\tilde V}(y)=\frac{1}{2\pi}
\int_{\R^2} e^{-ix\cdot y} \tilde V(x)\, dx.
}
Due to the radial symmetry of $\tilde V$, our transform gives radially symmetric function $\hat{\tilde V}(y)=\hat{\tilde V}(s)$, where $s=|y|$, that satisfies
\eq{\label{f1}
\hat{\tilde V}(s)&=\frac{1}{2\pi}\int_{0}^{2\pi}\int_0^\infty  e^{-isr \cos(\theta)}\tilde V(r) 
r\, dr\, d\theta\\
&=\int_0^R \tilde V(r) J_0(sr) r\, dr
=\frac{\kappa}{2}\int_0^{sR}\left[\lr{\frac{h}{s}-l_0}^2-(R-l_0)^2\right]J_0(h)\frac{h}{s^2}\, dh\\
&=
\frac{\kappa(2l_0-R)R}{2s^2}\int_{0}^{sR}h J_0(h)\, dh
-\frac{\kappa l_0}{s^3}\int_{0}^{sR}h^2 J_0(h)\, dh
+\frac{\kappa}{2 s^4}\int_{0}^{sR}h^3 J_0(h)\, dh,
}
where $J_0$ is the Bessel function of the first kind of order $0$.
In order to compute  integrals of the type  $\int_{0}^{H}h^\alpha J_0(h)\, dh$ for $\alpha=1,2,3$, we recall the Maclaurin series for the Bessel function of order~$i$
\eqh{
J_i(x)=\sum_{m=0}^\infty \frac{(-1)^m}{m! \Gamma(m+1+i)}\lr{\frac{x}{2}}^{2m+i};
}
then it is easy to check the following relations
\eqh{
x^{i+1} J_i(x)=(x^{i+1} J_{i+1}(x))', \quad x(J_{i-1}(x)+J_{i+1}(x))=2i J_i(x).
}
And so we easily compute
\eqh{
\int_{0}^{sR}h J_0(h)\, dh&=\int_{0}^{sR}\lr{h J_1(h)}'\, dh=sRJ_1(sR);\\
\int_{0}^{sR}h^3 J_0(h)\, dh&= (sR)^3 J_1(sR)-2\int_0^{sR} h^2J_1(h)\, dh\\
&= (sR)^3 J_1(sR)-2\int_0^{sR} \lr{h^2 J_2(h)}'\, dh\\
&= (sR)^3 J_1(sR)-2 (sR)^2 J_2(sR)=2(sR)^2 J_0(sR)+sR((sR)^2-4)J_1(sR).
}
However, the second integral on the r.h.s. of \eqref{f1} is more complicated, we have
\eqh{
\int_0^{sR} h^2 J_0(h)\, dx= (sR)^2J_1(sR) -\frac{\pi sR}{2}\lr{J_1(sR) H_0(sR)-J_0(sR)H_1(sR)},
}
where $H_0, H_1$ are the Struve functions defined by
\eqh{
H_i(x)=\sum_{m=0}^\infty\frac{(-1)^m}{\Gamma(m+3/2)\Gamma(m+i+3/2)}\lr{\frac{x}{2}}^{2m+i+1}.
}
Plugging the above formulas into \eqref{f1} we obtain
\eqh{
\hat{\tilde V}(s)&=\kappa\lr{J_0(sR)\frac{R^2}{s^2}-J_1(sR)\frac{2R}{s^3}+\frac{\pi R l_0}{2s^2}\left[J_1(sR)H_0(sR)-J_0(sR)H_1(sR)\right]}.
}
Therefore, the general  equation \eqref{tr_f} has now the following form
\eqh{
&\pt\log \hat f(t,y)\\
&\quad=-y^2
-2\pi f_\star\frac{\nu_f}{\nu_d}
\lr{J_0(|y|R){R^2}-J_1(|y|R)\frac{2R}{|y|}+\frac{\pi R l_0}{2}\left[J_1(|y|R)H_0(|y|R)-J_0(|y|R)H_1(|y|R)\right]}.
}
We now write an explicit form of the solution emanating from the initial data $f(0)=f_0$
\eqh{\hat f (t,y)=\hat f_0(y)e^{-G(y) t},}
where  the exponent $G= G(y,R,l_0, \kappa,\nu_f,\nu_d, f_\star)$ is given by
\eq{\label{Gy}
&G=y^2\\
&+2\pi f_\star\frac{\nu_f}{\nu_d} 
\lr{J_0(|y|R){R^2}-J_1(|y|R)\frac{2R}{|y|}+\frac{\pi R l_0}{2}\left[J_1(|y|R)H_0(|y|R)-J_0(|y|R)H_1(|y|R)\right]}.}
From Lemma \ref{L:un} we know exactly when $G$ ceases to be nonnegative close to $y=0$. Let us now see what happens slightly further from the origin. To this purpose, we rewrite \eqref{Gy} in the following form
\eqh{
G(z) R^2= z^2+2\pi f_\star\frac{\nu_f}{\nu_d}
{R^4}\lr{\frac{\pi l_0}{2R}\left[J_1(z)H_0(z)-J_0(z)H_1(z)\right]-J_2(z)},}
where we denoted $z=|y| R$.
To investigate the minima of $G(z)$ we check the minima of another function, namely
\eq{\label{defF}
F^{\alpha,\beta}(z)=G(z)R^2=z^2+\beta\lr{\frac{\pi \alpha}{2}\left[J_1(z)H_0(z)-J_0(z)H_1(z)\right]-J_2(z)},
}
where the parameters $\alpha, \beta>0$ are related to $R,l_0,\kappa,\nu_f,\nu_d, f_\star$ in the following way
\eq{\label{ab}
 \alpha=\frac{l_0}{R}, \quad \beta=\frac{2\pi\kappa f_\star\nu_f R^4}{ \nu_d}.
}
The interesting range for parameter $\alpha$ is $[0,1]$ and for the parameter $\beta$ we take $[0,\infty)$.
Below we present the graphs of the two functions $\frac{\pi}{2}\left[J_1(z)H_0(z)-J_0(z)H_1(z)\right]$ and $-J_2(z)$ that are included in the definition of $F^{\alpha,\beta}(z)$ from \eqref{defF}.
\begin{center}
\includegraphics[scale=0.4]{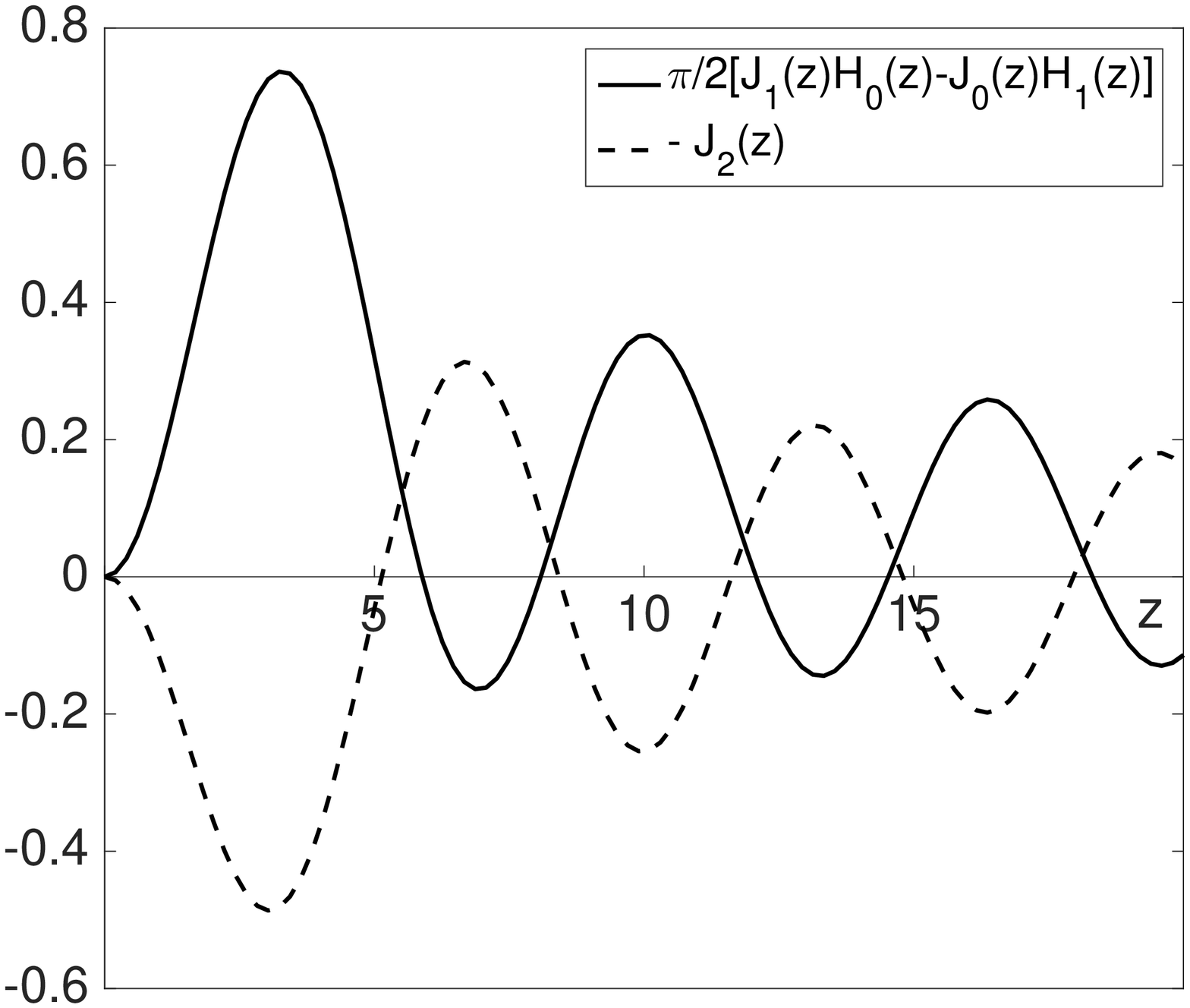}
\end{center}
Note, that from \eqref{defF} it is clear that $F^{\alpha,\beta}(0)=0$ for all values of $\alpha,\beta$. On the other hand, the picture above suggests that changing the values of parameters $\alpha,\ \beta$ may cause that  $F^{\alpha,\beta}$ will achieve negative values. In particular, by choosing a sufficiently small value for parameter $\alpha$
we would get a negative value of $\frac{\pi \alpha}{2}\left[J_1(z)H_0(z)-J_0(z)H_1(z)\right]-J_2(z)$ close to $z=0$.  This is nothing else than rephrasing the criterion from Lemma \ref{L:un} in terms of $\alpha$ and $\beta$.
\begin{prop}
Let $\alpha$ and $\beta$ be given as in \eqref{ab}, then if $(\alpha,\beta)\in U_{\R^2}$, where
\eqh{ U_{\R^2}=\left\{(\alpha,\beta)\in [0,1]\times [0,\infty): \alpha<\frac{3}{4},\  \beta>\frac{24}{3-4\alpha}\right\},}
then the steady state $f_\star$ is unstable, otherwise it is stable.
\end{prop}
\pf Instability of the steady state follows as previously from expansion of $F^{\alpha,\beta}(z)$ in the neighbourhood of $z=0$. After a bit lengthy but straightforward calculations we obtain
\eqh{
F^{\alpha,\beta}(z)= \lr{4+\beta\frac{2\alpha}{3}-\beta\frac{1}{2}}\lr{\frac{z}{2}}^2+O(z^4).
}
Finally, we see that taking $\alpha<\frac{3}{4}$ we can always find sufficiently large $\beta$  (i.e. $\beta>\frac{24}{3-4\alpha}$), so that the first term is negative and hence, for small enough $z$ the whole $F^{\alpha,\beta}(z)$ is negative as well.
The fact that for parameters $(\alpha,\beta)\notin U_{\R^2}$, the steady state is stable is shown numerically. On the picture below, we present the minimum of $F^{\alpha,\beta}$ with respect to $z$, i.e.
$$F^{\alpha,\beta}_{min}=\min_{z\in[0,10]}F^{\alpha,\beta}(z)$$
as a function of parameters  $\alpha,\beta$.
\begin{center}
\includegraphics[scale=0.35]{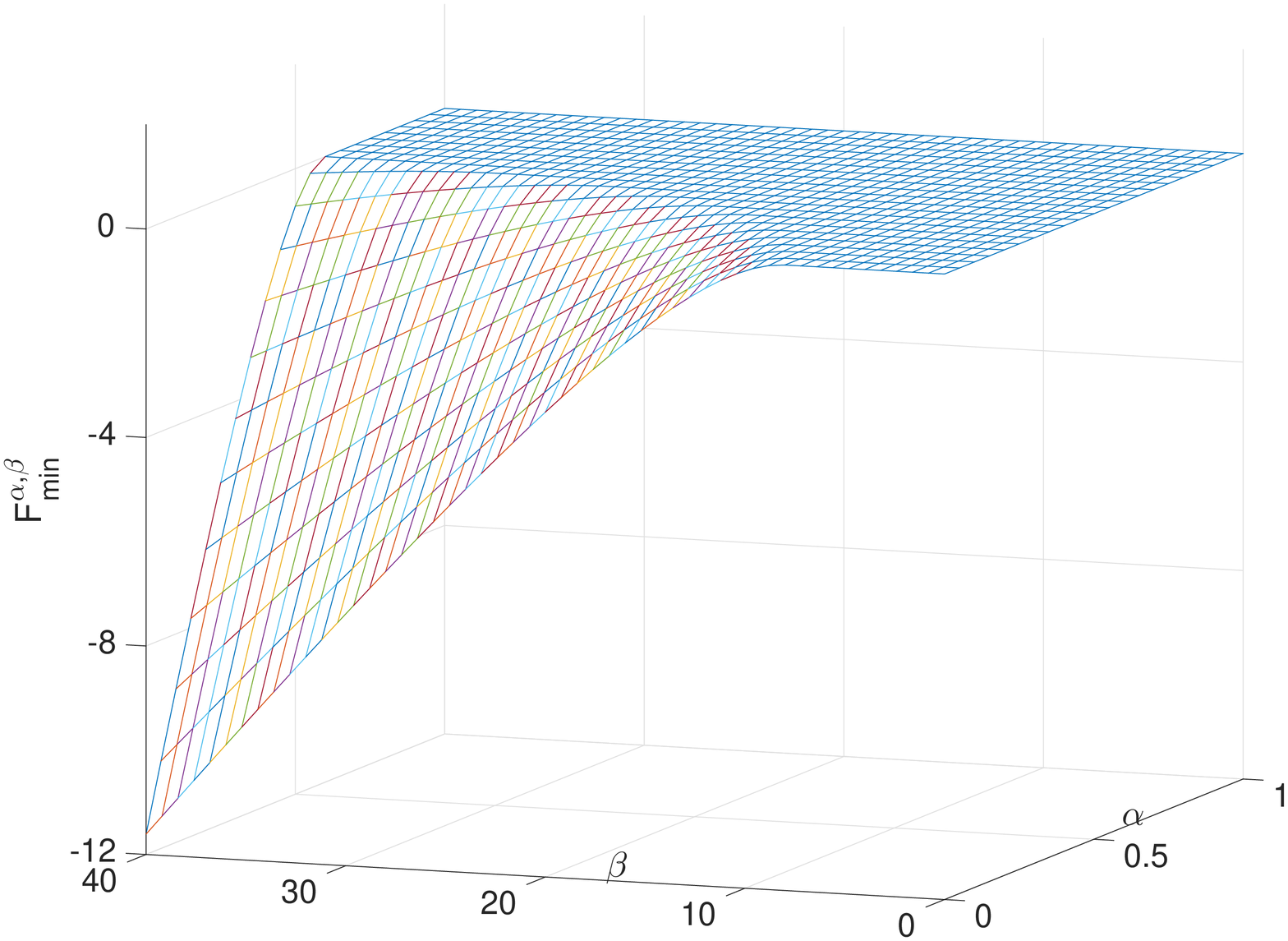}
\end{center}
The flat region corresponds to the parameter configuration that causes that the minimum of $F^{\alpha,\beta}(z)$ is attained at $z=0$ and is equal to $0$. $\Box$


\subsection{Linear analysis in the spacially-periodic case}\label{Sec:per}

Let us now investigate the same equation \eqref{main_per} but in the case of the space periodic domain. We will check an influence of  the size of the domain on the stability of stationary solutions. The analysis of  what happens with the solution in the unstable regime, but close to the instability  threshold will be presented in the next section.

We start by expanding our solution $f(x)$, for $x=(x_1,x_2)\in[-L_1,L_1]\times[-L_2,L_2]$ into the Fourier series. Introducing the shorthand notation for the Fourier modes
\eq{ e_{k_1,k_2}=\exp{\left[i\pi \lr{\frac{k_1x_1}{L_1}+\frac{k_2x_2}{L_2}}\right]},\label{not_e}} 
we may write
\eqh{
f(x_1,x_2)=\sum_{k_1,k_2\in \mathbb{Z}}\hat f_{k_1,k_2}e_{k_1,k_2},
}
where the Fourier coefficients $\hat f_{k_1,k_2}$ are given by
\eqh{
\hat f_{k_1,k_2}=\frac{1}{4L_1 L_2}\int_{-L_2}^{L_2}\int_{-L_1}^{L_1} f(x_1,x_2)e_{-k_1,-k_2}\, dx_1\, dx_2.}
Recall that we have the following properties for the Fourier coefficients of the derivatives of functions
\eqh{
\widehat{\partial^n_{x_1}f}_{k_1,k_2}=\lr{-i\frac{\pi k_1}{L_1}}^n\hat f_{k_1,k_2},\quad \widehat{\partial^n_{x_2}f}_{k_1,k_2}=\lr{-i\frac{\pi k_2}{L_2}}^n\hat f_{k_1,k_2}
}
and of the convolution of functions
\eqh{
\widehat{f\ast g}_{k_1,k_2}=\widehat{\left[\int_{-L_2}^{L_2}\int_{-L_1}^{L_1} f({x}
-{y})g({y})\, d{y}\right]}_{k_1,k_2}=4L_1L_2 \hat f_{k_1,k_2}\hat g_{k_1,k_2}.
}
Therefore, multiplying both sides of linearized system \eqref{main_per} by $\frac{1}{4L_1L_2} e_{-k_1,-k_2}$ and integrating over $[-L_1,L_1]\times[-L_2,L_2]$, we obtain
\eq{\label{2Df}
\pt\hat f_{k_1,k_2}=-{\pi^2}\lr{\frac{k_1^2}{L_1^2}+\frac{k_2^2}{L_2^2}} \hat f_{k_1,k_2}
-f_\star\frac{\nu_f}{\nu_d}{\pi^2}\lr{\frac{k_1^2}{L_1^2}+\frac{k_2^2}{L_2^2}} 4L_1L_2\hat{\tilde V}_{k_1,k_2}\hat f_{k_1,k_2}.
}
This time $f_\star$ can be interpreted as a probability measure, thus from now on, we will take $f_\star=\frac{1}{4L_1L_2}$ that on the rectangle $[-L_1,L_1]\times[-L_2,L_2]$ integrates to one, and so, for any $k_1,k_2\in \mathbb{Z}$, we obtain
\eqh{
\hat f_{k_1,k_2}(t)=\hat f_0(k_1,k_2)e^{-G_{k_1,k_2}t},
}
where
\eqh{
G_{k_1,k_2}={\pi^2}\lr{\frac{k_1^2}{L_1^2}+\frac{k_2^2}{L_2^2}}+ \frac{\nu_f}{\nu_d}{\pi^2}\lr{\frac{k_1^2}{L_1^2}+\frac{k_2^2}{L_2^2}}\hat{\tilde V}_{k_1,k_2}.
}
To compute $\hat{ \tilde V}_{k_1,k_2}$ in the case when $R<\min\{L_1,L_2\}$ we write
\eqh{
\hat{\tilde V}_{k_1,k_2}&=\frac{1}{4L_1L_2}\int_{0}^{2\pi}\int_0^R  e^{-i\pi\sqrt{{\frac{k_1^2}{L_1^2}+\frac{k_2^2}{L_2^2}}}r \cos\theta}\tilde V(r) 
r\, dr\, d\theta\\
&=\frac{\pi}{2L_1L_2}\int_0^R \tilde V(r) J_0\lr{\pi\sqrt{{\frac{k_1^2}{L_1^2}+\frac{k_2^2}{L_2^2}}}r} r\, dr}
and the last integral can be computed exactly as in the previous section so that we get
\eq{\label{fourV}
\hat{\tilde V}_{k_1,k_2}&=\frac{\kappa\pi}{2L_1L_2}\lr{\frac{\pi R^3 l_0}{2z_{k_1,k_2}^2}\left[J_1(z_{k_1,k_2})H_0(z_{k_1,k_2})-J_0(z_{k_1,k_2})H_1(z_{k_1,k_2})\right]-J_2(z_{k_1,k_2})\frac{R^4}{z_{k_1,k_2}^2}},
}
where we denoted 
\eq{\label{zjk}
z_{k_1,k_2}=\pi R\sqrt{{\frac{k_1^2}{L_1^2}+\frac{k_2^2}{L_2^2}}},}
and so
\eqh{
F^{\alpha,\beta}(z_{k_1,k_2})&=G_{k_1,k_2}R^2\\
&=z_{k_1,k_2}^2+\beta\lr{\frac{\pi\alpha}{2}\left[J_1(z_{k_1,k_2})H_0(z_{k_1,k_2})-J_0(z_{k_1,k_2})H_1(z_{k_1,k_2})\right]-J_2(z_{k_1,k_2})},}
for parameters $\alpha$ and $\beta$ such that
\eqh{
 \alpha=\frac{l_0}{R}, \quad \beta=\frac{\pi\kappa \nu_f R^4}{ 2\nu_dL_1L_2}.
}
Note that these are the same parameters as in \eqref{ab} with $f_\star=\frac{1}{4L_1L_2}$. Moreover, function $F^{\alpha,\beta}$ has the same form as in the whole space case  \eqref{defF}, but is evaluated only at the discrete set of points $z_{k_1,k_2}$ $k_1,k_2\in\mathbb{Z}$.
We know already that for continuous arguments $z\in[0,\infty)$ there is a phase transition curve $\beta(\alpha)=\frac{24}{3-4\alpha}$. The proof of this fact was based on finding a negative value of $F^{\alpha,\beta}(z)$ sufficiently close to $z=0$. Here,  however, the discrete variable $z_{k_1,k_2}$ depends on the size of the domain and it may happen that $F^{\alpha,\beta}(z_{k_1,k_2})$ for all $k_1,k_2\in \mathbb{Z}$ is always positive even if $F^{\alpha,\beta}(z)$ does attain negative value. Indeed, we have the following proposition.
\begin{prop}\label{Prop:abL}
For a nonempty subset of parameters  $(\alpha, \beta)\in U_{\R^2}$, there exist $L_1,L_2\in[R,\infty)$, such that 
$f_\star=\frac{1}{4L_1L_2}$ 
is a stable solution of \eqref{main_sys1}.
\end{prop}

\pf The proof of this fact is again numerical. The following graphs illustrate the function $F^{\alpha\beta}(z)$ in the unstable range of $\alpha,\beta$ ($\alpha=0.5, \beta=25)$  and the zoom of the graph in the neighbourhood of $0$.
\begin{center}
\includegraphics[scale=0.35]{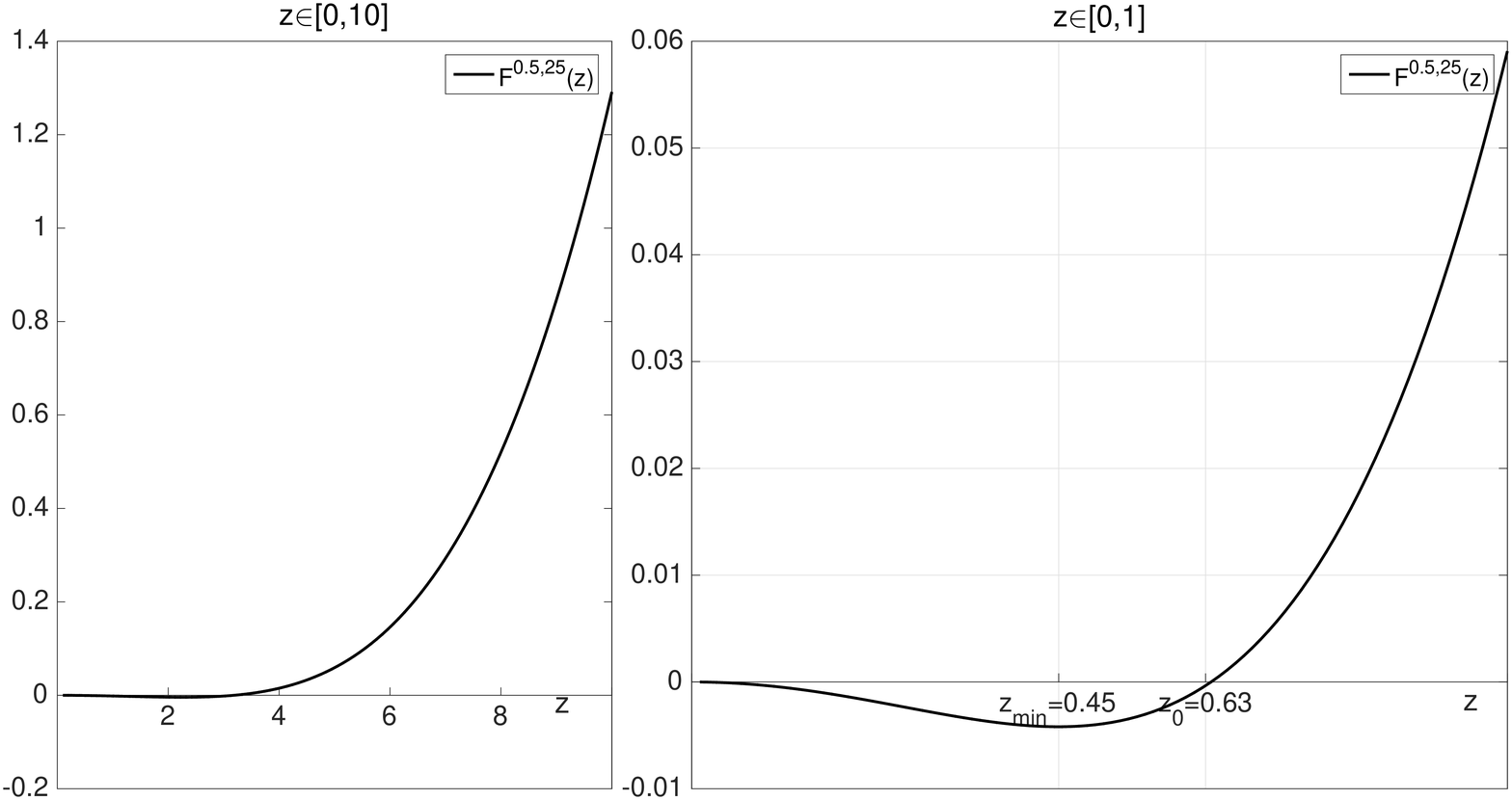}
\end{center}
Note, that if $z_{1,0}, z_{0,1}$ are larger than $z_0=0,63$ the steady state $f_\star$ will not be affected by the unsteady modes.  Below we present the positions of minima of function $F^{\alpha,\beta}(z)$, $z_{min}(\alpha,\beta)$
and the positions of zero of $F^{\alpha,\beta}(z)$, $z_{0}(\alpha,\beta)$. 

\vspace{-0.3cm}

\begin{center}
\includegraphics[scale=0.5]{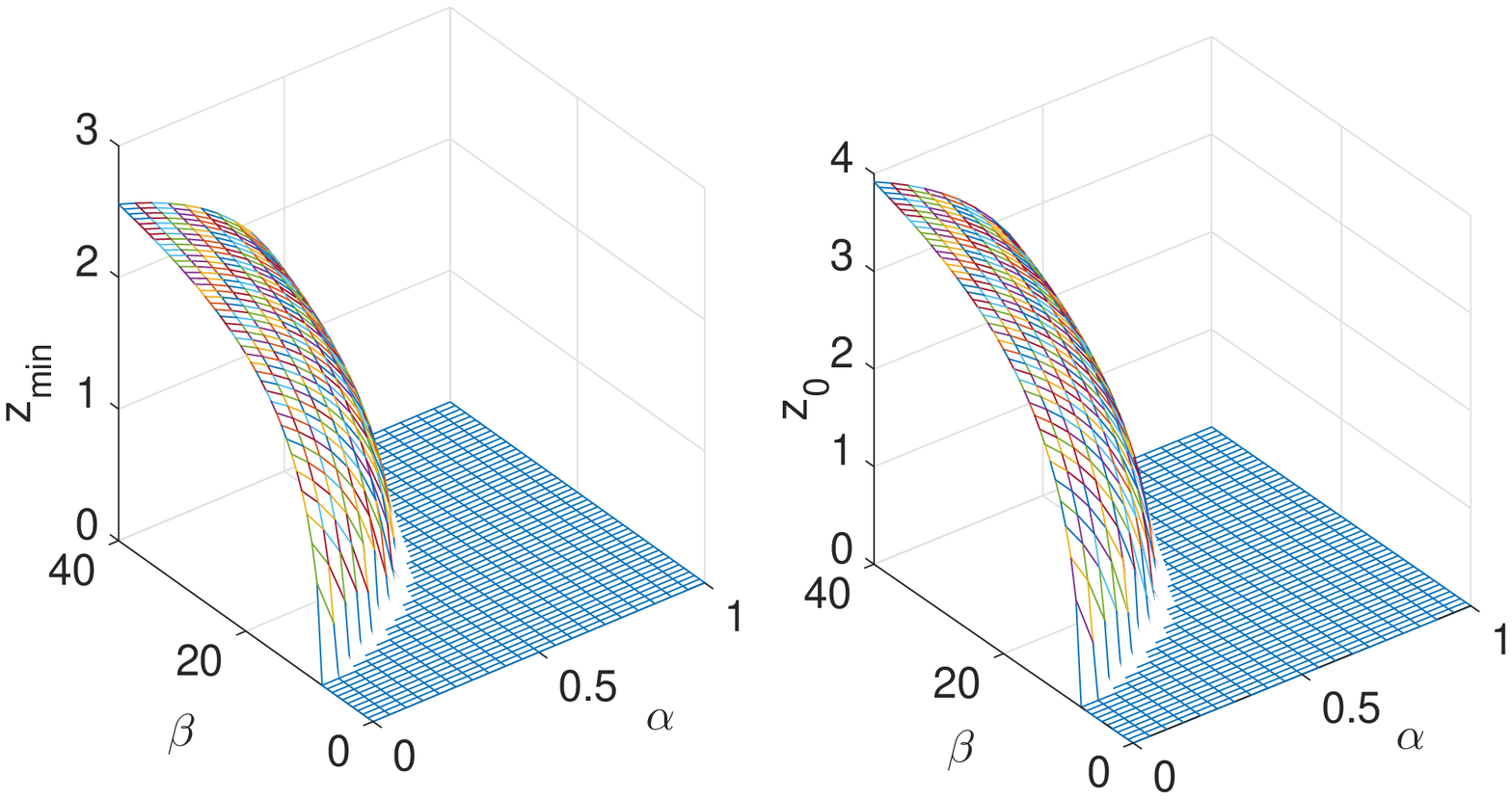}
\end{center}

\vspace{-1.5cm}

We see in particular, that  $z_0(\alpha,\cdot)$ is a monotonically increasing function, while $z_0(\cdot,\beta)$ is monotonically decreasing.
However, from \eqref{zjk} we get
$$z_{1,0}=\frac{\pi R}{L_1}\leq \pi\quad \mbox{and}\quad z_{0,1}=\frac{\pi R}{L_2}\leq \pi,$$
therefore, the statement can be fulfilled for example for $L_1=L_2=R$ and $\alpha^*,\beta^*$ such that
\eq{\label{c01}
z_0(\alpha^*,\beta^*)< \pi,
}
since $|z_0(\alpha,\beta)|\geq |z_{min}(\alpha,\beta)|$, the pair of parameters $(\alpha^*,\beta^*)\in U_{\R^2}$. $\Box$

The condition \eqref{c01} can be rephrased as follows
\eqh{
F^{\alpha,\beta}(\alpha^*,\beta^*)(\pi)>0,
}
which gives $\beta^*(0.7332\alpha^*-0.4854)>-9.8696$. This means in particular that   $\alpha^*\in (0.5499,0.75)$ and any $\beta^*\in[0,\infty)$ the stationary solution $f_\star=\frac{1}{4R^2}$ is a stable solution to \eqref{main_sys1} on a periodic box $[-R,R]^2$.

Using the same argument, we can also show the reverse statement to Proposition \ref{Prop:abL}, namely:
\begin{prop}\label{Prop:Lab}
For every $L_1,L_2\in[R,\infty)$, there exists a nonempty subset of parameters  $(\alpha, \beta)\in U_{\R^2}$, such that
$f_\star=\frac{1}{4L_1L_2}$ 
is a stable solution of \eqref{main_sys1}.
\end{prop}


\section{Nonlinear stability analysis of the steady-state}\label{sec:non}


\subsection{Preliminaries}

The purpose of this section is to investigate  the qualitative behavior of the model beyond the linear level.
We will choose the parameters $\alpha,\beta$ in the unstable regime, but close to the stability/instability threshold. 
In particular, the instability will be associated only with the first nontrivial modes, and the instability rate will be assumed small. As we saw in the previous section this can be guaranteed by the appropriate choice of the size of periodic domain.

The analysis will be made for periodic domains of two types: the rectangular periodic domain, and the square periodic domain. As we will see below, in the case when one side of the periodic domain is larger then the other, we may select only one unstable mode and reduce the analysis to a one-dimensional problem. For the case of a square box, 
the extra symmetry induces a degeneracy of the unstable mode. 
In both cases we give precise conditions for continuous and discontinuous phase transitions.
In the end of this section we also provide numerical verification of these conditions for the Hooke potential. However, we would like to emphasize that the theoretical results presented in this section are applicable to much wider class of potentials.
Our starting point is \eqref{main_sys1}, that we recall here for convenience 
\eq{\label{eq:main}
\partial_t f=\lap f+ \gamma\Grad \cdot (f \Grad (\tilde V\ast f)),}
with $\gamma=  \frac{\nu_f}{\nu_d}$.


\subsection{The rectangular case - non degenerate}

We start our analysis from the simpler case when the periodic domain is rectangular
$$(x_1,x_2)\in[-L_1,L_1]\times [-L_2,L_2],\quad \mbox{such that}\quad L_1>L_2,$$ 
and that only the modes $(\pm1,0)$ are unstable, all the others are stable.
Having in mind the argument from the previous section, this is possible for some $(\alpha^*,\beta^*)\in U_{\R^2}$ provided
$$z_{1,0}< z_0(\alpha^*,\beta^*)<z_{2,0},\quad \mbox{and}\quad z_{0,1}>z_0(\alpha^*,\beta^*).$$
Looking at the problem from the perspective of stable and unstable modes, we see that an analogous condition can be deduced directly from \eqref{2Df}. Namely, the eigenvalue associated with the first mode in the direction $x_1$ should be the only positive one. This results in the conditions:
 \begin{subequations}
\eq{
\lambda=\lambda_{\pm 1,0}=-\frac{\pi^2}{L_1^2}\lr{1+\gamma\hat{\tilde V}_{1,0}}>0, \label{lambda1}}
\eq{\lambda_{k_1,k_2}=-{\pi^2}\lr{\frac{k_1^2}{L_1^2}+\frac{k_2^2}{L_2^2}} \lr{1+\gamma\hat{\tilde V}_{k_1,k_2}}<0,\quad \mbox{for}\ (k_1,k_2)\neq(\pm1,0).}
 \end{subequations}
Recalling notation \eqref{not_e}, the unstable modes are then: 
$$e_{1,0}=e^{\frac{i\pi x_1}{L_1}}\quad \mbox{and} \quad e_{-1,0} = e^{\frac{-i\pi x_1}{L_1}}.$$
We now want to check what happens with the constant steady state after passing the instability threshold. We could, for example, think of fixing the parameter $\alpha$ according to Proposition \ref{Prop:Lab} and slowly increase parameter $\beta$ by changing the value of $R$. Alternatively, one can identify the instability threshold with changing the sign of $\lambda$ -- this is the standard strategy in bifurcation theory, and the one we follow here.

After crossing the instability threshold, one expects that the solution to the nonlinear problem behaves for short time like the linearized solution, that is, an exponential in time times the unstable mode:
$$f=f_\star+A(t)e_{1,0}+A^\ast(t) e_{1,0},\quad \mbox{with}\quad A(t)\propto e^{\lambda t}.$$
Then, if $A(t)$ remains small, one can hope to expand  the solution into power series of $A(t)$
$$f=f_\star+A(t)e_{1,0}+A^\ast(t) e_{1,0}+O(A(t)^2)$$
The goal is then to find a reduced equation for $A(t)$ that would allow us to understand the dynamics of the solution just by analyzing an ODE for $A(t)$ (central manifold reduction). The unstable eigenvalue is real, and the system is translation-symmetric, hence we expect a pitchfork bifurcation when $\lambda$ changes sign from "$-$" to "$+$", with two possible scenarios: \\
$\bullet$ A {\it supercritical bifurcation}: $A(t)$ first grows exponentially, but then $f$ tends to an almost homogeneous stationary state,\\
$\bullet$ A {\it subcritical bifurcation}: $A(t)$ grows exponentially until it leaves the perturbative regime, then the final state may be very far from the original homogeneous state.

Instead of adopting a dynamical approach as done here, bifurcations for systems such as \eqref{eq:main} can be studied from a "thermodynamical" point of view, i.e. by looking at the minimizers of \eqref{eq:entropy}. This has been done in particular in \cite{ChPa}. The {\it second order phase transition} in \cite{ChPa} corresponds to the supercritical scenario described above, while the {\it first order phase transition} corresponds to the subcritical scenario. However, one should note that the dynamical bifurcation point (where $\lambda$ changes sign) does not coincide with the first order phase transition parameters; the dynamical bifurcation would rather be called a {\it spinodal point}  in thermodynamics, the language of  \cite{ChPa}.

The main result of this section provides a criterion allowing to distinguish these two cases.
\begin{thm}\label{Th1D}
Assume that $\lambda>0$ and that $\lambda_{k_1,k_2}<0$ for any $(k_1,k_2)\neq(\pm1,0)$. Then, there are two possibilities:\\
$\bullet$  for $2\hat{\tilde V}_{2,0}-\hat{\tilde V}_{-1,0}>0$ the steady state exhibits a supercritical bifurcation,\\
$\bullet$  for $2\hat{\tilde V}_{2,0}-\hat{\tilde V}_{-1,0}<0$ the steady state exhibits a subcritical bifurcation.\\
\end{thm}

\pf We now want to investigate the evolution of the perturbation $g$ of the constant steady state $f_\star$. Hence, the solution to \eqref{main_sys1} has the form $f=f_\star+\eta$. We denote the operator associated with the linearized equation \eqref{main_per}  by ${\cal L}(f)$, more precisely
 \eqh{
\partial_t \eta(t,x)=\lap_x \eta(t,x)+ 
 \gamma f_\star\lap_x( (\tilde V\ast \eta)(t,x)):= {\cal L}(\eta),
}
Note that ${\cal L}(\eta)$ with periodic boundary conditions is a self adjoint operator.
Next, we also distinguish the nonlinear part of \eqref{main_sys1} and we denote it by ${\cal N}(\eta)$, this gives
\eq{\label{lpnl}
\partial_t \eta= \mathcal{L}(\eta) +\mathcal{N}(\eta),
}
where
\eqh{
\mathcal{N}(\eta) = \mathcal{Q}(\eta,\eta),\qquad \mathcal{Q}(\eta_1,\eta_2) = \gamma\nabla\cdot\left( \eta_1\nabla (\tilde V\ast \eta_2) \right) }
with $\gamma=\frac{\nu_f}{\nu_d}$. In what follows we will need to compute the action of $\mathcal{L}$ and $\mathcal{Q}$ on the Fourier basis. We have 
\begin{eqnarray}
\mathcal{L}(e_{k_1,k_2}) &=& \left[-\pi^2\left(\frac{k_1^2}{L_1^2}+\frac{k_2^2}{L_2^2}\right)(1+\gamma\hat{\tilde{V}}_{k_1,k_2}) \right] e_{k_1,k_2}=\lambda_{k_1,k_2} e_{k_1,k_2}  \label{eq:Lfour}\\
\mathcal{Q}(e_{k_1,k_2},e_{l_1,l_2}) &=& -4L_1L_2\gamma \pi^2 \hat{\tilde{V}}_{l_1,l_2} \left(\frac{l_1(k_1+l_1)}{L_1^2}+\frac{l_2(k_2+l_2)}{L_2^2}\right)e_{k_1+l_1,k_2+l_2}
\label{eq:Qfour}
\end{eqnarray}
As mentioned above, at a linear order, $\eta$ moves on a vector space spanned by $e_{1,0},e_{-1,0}$
\begin{equation*}
\eta(t,x) = A(t) e_{1,0} +A^\ast(t) e_{-1,0}.
\label{eq:expg}
\end{equation*}
Furthermore, if the equation were linear, solution emanating from any initial condition would be quickly attracted towards this vector space. This follows from the fact that all the other modes of motion are stable.
For the nonlinear system, we expect that ${\rm span}(e_{1,0},e_{-1,0})$ will be deformed into some manifold. This manifold is tangent to ${\rm span}(e_{1,0},e_{-1,0})$ close to $\eta=0$, and can be parametrized by the projection of $\eta$ on this space as follows
\begin{equation}
\eta(t,x) = A(t) e_{1,0}+A^\ast(t) e_{-1,0} +H[A,A^\ast](x),
\label{eq:expg}
\end{equation}
with $H$ such that
\eq{\label{condH}
 H[A,A^\ast] =O(A^2, AA^\ast,(A^\ast)^2) \quad\mbox{and}\quad \langle e_{1,0},H\rangle=\langle e_{-1,0},H\rangle=0.} 
 Furthermore, from translation invariance we can write, using Lemma \ref{Lem:H} (see below):
\begin{equation*}
H[A,A^\ast]=\sum_{k_1\geq 0}   A^{k_1}   h_{k_1,0} (\sigma) e_{k_1,0} + \sum_{k_1< 0}   (A^\ast)^{-k_1}   h_{k_1,0} (\sigma) e_{k_1,0},
\end{equation*}
where 
\eqh{\sigma=|A|^2\quad \mbox{and} \quad h_{k,0}=h_{k,0}^0+\sigma h_{k,0}^1+\ldots.
}
The conditions \eqref{condH} imply that $h_{1,0}=h_{-1,0}=0$. Moreover, $h^0_{k_1,0}=0$ for $k_1=0,\pm1$, otherwise $H[A,A^\ast]$ would contain zero and first order terms in $A, A^\ast$. Hence, at the leading order, only the modes $(\pm 2,0)$ remain; more precisely
\eq{
\label{eq:expHb}
H[A,A^\ast] = A^2 h_{2,0}^0 e_{2,0} + (A^{\ast})^2 h_{-2,0}^0 e_{-2,0} + O((A,A^\ast)^3).}
Then, plugging \eqref{eq:expg} and \eqref{eq:expHb} into the definitions of ${\cal L}(\eta)$ and  ${\cal N}(\eta)$ we obtain
\eq{\label{eq:dtg0}
 \mathcal{L}(\eta) = A {\cal L}(e_{1,0}) +A^\ast {\cal L}(e_{-1,0}) + A^2 h_{2,0}^0{\cal L}(e_{2,0})+ (A^\ast)^2 h_{-2,0}^0  {\cal L}(e_{-2,0})
+O((A,A^\ast)^3),
}
and
\eq{
\mathcal{N}(\eta)=&
A^2\mathcal{Q}(e_{1,0},e_{1,0})+(A^\ast)^2\mathcal{Q}(e_{-1,0},e_{-1,0})\\
&+ A^3 h_{2,0}^0\left[\mathcal{Q}(e_{1,0},e_{2,0}) +\mathcal{Q}(e_{2,0},e_{1,0})\right] \\
&+|A|^2A  h_{2,0}^0 \left[\mathcal{Q}(e_{-1,0},e_{2,0}) +\mathcal{Q}(e_{2,0},e_{-1,0})\right]\\
&
+(A^\ast)^3 h_{-2,0}^0\left[\mathcal{Q}(e_{-1,0},e_{-2,0}) +\mathcal{Q}(e_{-2,0},e_{-1,0})\right]\\
&+|A|^2A^\ast h_{-2,k}^0\left[\mathcal{Q}(e_{1,0},e_{-2,0}) +\mathcal{Q}(e_{-2,0},e_{1,0})\right]\\
&+O((A,A^\ast)^4).
\label{eq:dtg1}
}
Therefore, the full dynamics of $\eta$ can be obtained by substituting the above formulas for ${\cal L}(\eta)$ and  ${\cal N}(\eta)$ into \eqref{lpnl}. On the other hand, differentiating \eqref{eq:expg} with respect to time, and using \eqref{eq:expHb} we have
\eq{
\partial_t \eta &= \dot{A} e_{1,0}+\dot{A^\ast} e_{-1,0}+ 2A\dot{A} h_{2,0}^0 e_{2,0}+2 A^\ast\dot{A^\ast} h_{-2,0}^0 e_{-2,0}+\pt O((A,A^\ast)^3).
\label{eq:dtg2}
}
We now equate expressions $\partial_t\eta=\eqref{eq:dtg0}+\eqref{eq:dtg1}$ and \eqref{eq:dtg2}, and compare Fourier mode by Fourier mode, and order  in $A$ by order in $A$. We start with the mode $e_{1,0}$. 
Taking the scalar product of the right hand sides of \eqref{eq:dtg0} and \eqref{eq:dtg1} with $e_{1,0}$, we get 
\eqh{\langle  e_{1,0},\mathcal{L}(\eta)\rangle=A\langle  e_{1,0},\mathcal{L}(e_{1,0})\rangle,}
and
\eqh{\langle e_{1,0},  \mathcal{N}(\eta)\rangle = |A|^2 Ah_{2,0}^0\langle e_{1,0}, \mathcal{Q}(e_{-1,0},e_{2,0}) +\mathcal{Q}(e_{2,0},e_{-1,0})\rangle+O((A,A^\ast)^4),
}
where $\langle u,v\rangle=1/(4L_1L_2)\int_{-L_2}^{L_2}\int_{-L_1}^{L_1} u^* v\, dx_1\ dx_2$.
Comparing these expressions with the projection of \eqref{eq:dtg2} on $e_{1,0}$ we obtain
\eq{\label{eq1}
&\dot{A}\langle  e_{1,0},e_{1,0}\rangle \\
&\quad= A\langle  e_{1,0},\mathcal{L}(e_{1,0})\rangle+|A|^2 Ah_{2,0}^0\langle e_{1,0}, \mathcal{Q}(e_{-1,0},e_{2,0}) +\mathcal{Q}(e_{2,0},e_{-1,0})\rangle +O((A,A^\ast)^4).}
So, using \eqref{eq:Lfour} and \eqref{eq:Qfour} we obtain
\eqh{
\dot{A} = A\lambda+|A|^2 Ah_{2,0}^0\gamma\pi^2\frac{4L_2}{L_1}\lr{\hat{\tilde V}_{-1,0}-2\hat{\tilde V}_{2,0}} +O((A,A^\ast)^4).}
The terms of the leading order in $A$ yield the linearized dynamics. 
To investigate the behaviour of $A$ at the non-linear level we need first to compute $h_{2,0}^0$: we do this by equating the Fourier coefficient $(2,0)$ in
$\partial_t\eta=\eqref{eq:dtg0}+\eqref{eq:dtg1}$  and \eqref{eq:dtg2}; we obtain
\[2A\dot{A} h_{2,0}^0 e_{2,0}=A^2 \lambda_{2,0}h_{2,0}^0e_{2,0}+A^2\mathcal{Q}(e_{1,0},e_{1,0})+O((A,A^\ast)^4),
\]
so, using \eqref{eq:Lfour} and \eqref{eq:Qfour} together with linear equation for  $A$, i.e. $\dot{A}=\lambda A$, we obtain
\[
2\lambda h_{2,0}^0=-\frac{4\pi^2}{L_1^2}\left(1+\gamma \hat{\tilde{V}}_{2,0}\right)  h_{2,0}^0-\frac{8\pi^2L_2}{L_1}\gamma\hat{\tilde{V}}_{1,0},
\]
and finally, since $\lambda_{2,0}<0$, for $\lambda\to0^+$ we formally get
\[
h_{2,0}^0 = -\frac{-2L_1L_2\gamma\hat{\tilde{V}}_{1,0}}{1+\gamma\hat{\tilde{V}}_{2,0}}.
\]
The reduced equation for $A$ \eqref{eq1} then reads:
\begin{equation}
\dot{A}=\lambda A +8\gamma^2\pi^2L_2^2\frac{\hat{\tilde{V}}_{1,0}}{1+\gamma\hat{\tilde{V}}_{2,0}}\left(2\hat{\tilde{V}}_{2,0}-\hat{\tilde{V}}_{-1,0} \right)|A|^2A
\label{eq:reduced}
\end{equation}

From the assumptions of Theorem \ref{Th1D} and \eqref{lambda1} it follows that $\hat{\tilde V}_{1,0}$ is negative, so if $2\hat{\tilde V}_{2,0}-\hat{\tilde V}_{-1,0}>0$ the coefficient in front of the third order term is negative. This means that $A(t)$ first grows exponentially, but then it saturates when the r.h.s. of \eqref{eq:reduced} is equal to zero. This happens for
\[
|A|=\frac{\sqrt{\lambda}}{2\sqrt{2}\gamma\pi L_2} \sqrt{\frac{1+\gamma\hat{\tilde{V}}_{2,0}}{|\hat{\tilde{V}}_{1,0}|(2\hat{\tilde{V}}_{2,0}-\hat{\tilde{V}}_{-1,0})}}.
\]
Therefore, if the last factor is bounded $|A|$ is of order $\sqrt{\lambda}$, so, taking $\lambda$ sufficiently small we assure that $A(t)$ remains small at the level of saturation, which justifies the validity of expansion \eqref{eq:expg}.

When $2\hat{\tilde V}_{2,0}-\hat{\tilde V}_{-1,0}<0$ the term of order $A^3$ does not bring any saturation. The growth thus goes on until $A(t)$ leaves the perturbative regime, and at this point the approach breaks down. 

This yields the hypothesis of Theorem \ref{Th1D}.  In order to conclude, we still need to justify that the manifold $H$ can be represented by \eqref{eq:expHb}, we will prove the following lemma. 
\begin{lemma}\label{Lem:H}
Let $H=H[A,A^\ast](x)$ be as specified above in \eqref{eq:expg}, then $\hat H_{0,0}[A,A^\ast]=0,\ \hat H_{\pm1,0}[A,A^\ast]=0$ and the other
Fourier coefficients of $H$ are of the form
\[
\hat{H}_{k_1,k_2}[A,A^\ast] =\left\{ 
\begin{array}{lll}
A^{k_1} h_{k_1,0}(\sigma)&{\mbox for}&k_1\geq 0,\ k_2=0\\
(A^\ast)^{-k_1} h_{k_1,0}(\sigma)&{\mbox for}&k_1<0,\ k_2=0 \\
0 &{\mbox for}& k_2 \neq 0,
\end{array}\right.
\]
for some unknown functions $h_{k_1,0}=h_{k_1,0}(\sigma)$, with $\sigma =AA^\ast$.
\end{lemma}
\pf From the definition $\hat{H}_{0,0}=0$, and $\hat{H}_{\pm 1,0}=0$ since $\langle e_{1,0},H\rangle=\langle e_{-1,0},H\rangle=0$.\\
Next, equation  \eqref{lpnl}  as well as the unstable manifold are invariant under translation $\tau_{x^0} :x\to x+x^0$ that act on functions as
\[
(\tau_{x^0}\cdot f)(x) = f(x-x^0),
\]
where $x=(x_1,x_2)$, $x^0=(x^0_1,x^0_2)$.
Therefore, for any $A$, there exists $\tilde A$ such that
\begin{equation*}
\tau_{x^0}\cdot(A e_{1,0}+A^\ast e_{-1,0}+H[A,A^\ast]) =\tilde A e_{1,0}+\tilde A^\ast e_{-1,0}+H[\tilde A,\tilde A^\ast] ,
\end{equation*}
meaning that
\eqh{
&A e^{-i\pi \frac{x^0_1}{L_1}} e_{1,0}+A^\ast  e^{i\pi\frac{x^0_1}{L_1}} e_{-1,0}+H[A,A^\ast](x-x_0)\\
 &\qquad=\tilde A e_{1,0}+\tilde A^\ast e_{-1,0}+H[\tilde A,\tilde A^\ast](x).
}
comparing the terms with $e_{1,0}$ we conclude that $\tilde A=A e^{-i\pi \frac{x^0_1}{L_1}}$ and subsequently 
$$H\left[A e^{-i\pi \frac{x^0_1}{L_1}},A^\ast e^{i\pi \frac{x^0_1}{L_1}}\right](x)=H[A,A^\ast](x-x_0).$$ 
In terms of Fourier coefficients, the last equality reads
\begin{equation}
\hat{H}_{k_1,k_2}\left[A e^{-i\pi \frac{x^0_1}{L_1}},A^\ast e^{i\pi \frac{x^0_1}{L_1}}\right]= e^{-i\pi \lr{\frac{k_1 x^0_1}{L_1}+\frac{k_2 x^0_2}{L_2}}} \hat{H}_{k_1,k_2}[A,A^\ast].
\label{eq:fourierHk}
\end{equation} 
Let us now expand $\hat{H}_{k_1,k_2}$ in a Taylor series: $\hat{H}_{k_1,k_2}[z,z^\ast]=\sum_{l_1,l_2\geq 0} c_{l_1,l_2} z^{l_1}(z^\ast)^{l_2}$, then \eqref{eq:fourierHk} reads
\eqh{
\sum_{l_1,l_2\geq0} c_{l_1,l_2} A^{l_1}(A^\ast)^{l_2}e^{-i\pi\frac{x^0_1}{L_1}(l_1-l_2)} = e^{-i\pi \lr{\frac{k_1 x^0_1}{L_1}+\frac{k_2 x^0_2}{L_2}}}\sum_{l_1,l_2\geq0} c_{l_1,l_2} A^{l_1}(A^\ast)^{l_2}.
}
The uniqueness of the expansion implies that $c_{l_1,l_2}=0$ unless $l_1-l_2=k_1,\ k_2=0$. Thus
\eqh{
\hat{H}_{k_1,0}[A,A^\ast] = A^{k_1} \sum_{l_2\geq 0} c_{k_1+l_2,l_2} |A|^{2l_2}.}
$\Box$

\medskip

This finishes the proof of Theorem \ref{Th1D}. $\Box$


\subsection{The square case - degenerate eigenvalues}

In this section we study a particular case of domain -- a periodic box, thus $L_1=L_2=L$. For simplicity, we take $L=\frac{1}{2}$. Again, the result is much more general and might be applied to much wider class of functionals than the Hooke potential from Section \ref{Sec:per}, provided one can select finitely many unstable modes. Here, due to the square symmetry, and assuming that the potential is isotropic, there will generically be one unstable mode in each direction denoted by
$$e_{1,0}=e^{2i\pi x_1}\quad\mbox{and}\quad  e_{0,1}=e^{2i\pi x_2},$$
together with their conjugates, associated with the same eigenvalue
$$\lambda= -4 \pi^2 \lr{1+\gamma\hat{\tilde V}_{1,0}}.$$
Our results in this case can be summarized as follows.
\begin{thm}\label{Th2D}
Assume that $\lambda>0$ and that $1+\gamma\hat{\tilde V}_{k_1,k_2}>0$  for any $k_1,k_2$ such that $|k_1|+|k_2|>1$.
Then, for 
\eq{\label{super2D}
 \frac{\hat{\tilde V}_{1,0}(2\hat{\tilde V}_{2,0}-\hat{\tilde V}_{-1,0})}{1+\gamma\hat{\tilde V}_{2,0}}<-\left| 2 \frac{\hat{\tilde V}_{1,0}\hat{\tilde V}_{1,1}}{1+\gamma\hat{\tilde V}_{1,1}}\right|
}
the steady state exhibits a supercritical bifurcation. If the inequality is opposite, the steady state exhibits a subcritical bifurcation.
\end{thm}

\pf Following the same strategy as for the 1D case we expand the perturbation $\eta$ on the unstable manifold:
\begin{equation*}
\eta(t,x,y) = A(t)e_{1,0} +A^\ast(t) e_{-1,0} +B(t)e_{0,1} +B^\ast(t) e_{0,-1} +H[A,A^\ast,B,B^\ast](x,y),
\end{equation*}
therefore
\eqh{
\pt\eta(t,x,y)= \dot{A}e_{1,0} +\dot{A^\ast} e_{-1,0} +\dot{B}e_{0,1} +\dot{B^\ast} e_{0,-1} +\pt H[A,A^\ast,B,B^\ast](x,y).
}
Alike in Lemma \ref{Lem:H}, we can deduce that $H$ has the following structure
\eq{
\label{eq:expH2D}
H =& A^2 h_{2,0} e_{2,0}+ (A^{\ast})^2 h_{-2,0} e_{-2,0} +B^2 h_{0,2} e_{0,2}+ (B^{\ast})^2 h_{0,-2} e_{0,-2}  \\
& +AB h_{1,1} e_{1,1}+A^\ast B h_{-1,1} e_{-1,1} +AB^\ast h_{1,-1} e_{1,-1}+A^\ast B^\ast h_{-1,-1} e_{-1,-1} \\
&+O((A,A^\ast,B,B^\ast)^3).
}
We compute now the non linear term $\mathcal{N}(\eta)$ at order $A^2,B^2$ (we use here the properties of $\tilde{V}$: $\hat{\tilde{V}}_{k_1,k_2}=\hat{\tilde{V}}_{k_1,-k_2}=\hat{\tilde{V}}_{-k_1,k_2}=\hat{\tilde{V}}_{k_2,k_1}$):
\eqh{
\mathcal{N}(\eta) = -8\gamma \pi^2 \hat{\tilde V}_{1,0}\left[A^2  e_{2,0}+ B^2 e_{0,2} +AB e_{1,1} +A^\ast B e_{1,-1}
 +{\rm c.c.}\right] +O\left((A,A^\ast,B, B^\ast)^3\right)
}
The procedure is the same as before. The leading order for the dynamics of $A,B$ is the linear evolution:
\[
\dot{A} =\lambda A +O((A,B)^3),\qquad \dot{B} =\lambda B +O((A,B)^3).
\]
We expand in powers of $\sigma_A=|A|^2, \sigma_B=|B|^2$ the $h_{kl}$ coefficients that appear in \eqref{eq:expH2D}, and keep only the leading order $h_{k,l}^0$, which are some constants to be computed. From comparison of $(2,0)$, $(1,1)$ and $(1,-1)$ modes respectively, order $(A,B)^2$ yields the equations for $h_{\pm 2,0}^0, h_{0,\pm 2}^0, h_{\pm 1,\pm 1}^0$:
\eqh{
(2\lambda -\lambda_{2,0})h^0_{2,0} &= -8\gamma\pi^2 \hat{\tilde V}_{1,0},\\
(2\lambda -\lambda_{1,1})h^0_{1,1} &= -8\gamma\pi^2 \hat{\tilde V}_{1,0},\\
(2\lambda -\lambda_{1,-1})h^0_{1,-1} &=  -8\gamma\pi^2 \hat{\tilde V}_{1,0}. 
}
Solving the above equations, and letting $\lambda\to0$, we obtain:
\eqh{
h_{2,0}^0 = -\frac{\gamma \hat{\tilde V}_{1,0}}{2(1+\gamma\hat{\tilde V}_{2,0})}, \qquad
h_{1,1}^0 = -\frac{\gamma\hat{\tilde V}_{1,0}}{1+\gamma\hat{\tilde V}_{1,1}}, \qquad
h_{1,-1}^0 = -\frac{\gamma\hat{\tilde V}_{1,0}}{1+\gamma\hat{\tilde V}_{1,1}}.
}
The other relevant $h_{i,j}^0$ coefficients in \eqref{eq:expH2D} are obtained by complex conjugation.
Finally, including the terms of order $(A,B)^3$ for the Fourier modes $(1,0)$ and $(0,1)$ we obtain the sought reduced equations for evolution of $A$ and $B$, namely
\eqh{
\left\{\begin{array}{rl}
\dot{A} = &\lambda A +|A|^2 Ah_{2,0}^0\left[\mathcal{Q}(e_{-1,0}, e_{2,0})+\mathcal{Q}(e_{2,0},e_{-1,0}) \right]\\
&+|B|^2 Ah_{1,-1}^0\left[\mathcal{Q}(e_{0,1}, e_{1,-1})+\mathcal{Q}(e_{1,-1},e_{0,1}) \right]\\
&+|B|^2 Ah_{1,1}^0\left[\mathcal{Q}(e_{1,1}, e_{0,-1})+\mathcal{Q}(e_{0,-1},e_{1,1}) \right]+O((A,A^\ast, B, B^\ast)^4), \\
\dot{B} = &\lambda B +|B|^2Bh_{0,2}^0\left[\mathcal{Q}(e_{0,-1}, e_{0,2})+\mathcal{Q}(e_{0,2},e_{0,-1}) \right]\\
&+|A|^2 Bh_{-1,1}^0\left[\mathcal{Q}(e_{1,0}, e_{-1,1})+\mathcal{Q}(e_{-1,1},e_{1,0}) \right]\\
&+|A|^2 Bh_{1,1}^0\left[\mathcal{Q}(e_{-1,0}, e_{1,1})+\mathcal{Q}(e_{1,1},e_{-1,0}) \right]+O((A,A^\ast, B, B^\ast)^4),
\end{array}
\right.
}
or equivalently
\eq{\label{sysAB}
\left\{\begin{array}{l}
\dot{A} = \lambda A +c|A|^2 A+d|B|^2A, \\
\dot{B} = \lambda B +c|B|^2 B+d|A|^2B, 
\end{array}
\right.
}
where we denoted 
\eq{\label{cd}
c =2\gamma^2\pi^2 \frac{\hat{\tilde V}_{1,0}(2\hat{\tilde V}_{2,0}-\hat{\tilde V}_{1,0})}{1+\gamma\hat{\tilde V}_{2,0}},\qquad d= 8\gamma^2\pi^2 \frac{\hat{\tilde V}_{1,0}\hat{\tilde V}_{1,1}}{1+\gamma\hat{\tilde V}_{1,1}}~.
}
The analysis of the two-dimensional system requires slightly more effort than the analysis of the one-dimensional case from the previous section. The steady states of the system \eqref{sysAB} are determined by
\eqh{
 (\lambda +c|A|^2+d|B|^2)A=0,\quad \mbox{and}\quad (\lambda +c|B|^2+d|A|^2)B=0.
}
If $c<0$, there are steady states with $A=0$ or $B=0$; it is easy to see that they are unstable. If $c+d<0$, there are other steady states, with $A=A_{\rm st}\neq 0$ and $B=B_{\rm st}\neq 0$. The modulus of $A_{\rm st}$ and $B_{\rm st}$ is fixed, but their phase is undetermined:
\[
|A_{\rm st}|=|B_{\rm st}|=\sqrt{\frac{\lambda}{-(c+d)}}.
\]
In order to check stability of the above steady states, we investigate the linearization of system \eqref{sysAB}, around $(A_{\rm st}, B_{\rm st})$. We take for simplicity 
$A_{\rm st}$ and $B_{\rm st}$ real in the following; by translation symmetry the result does not depend on the phases we choose. Furthermore, one checks easily that
the linearized equations for the imaginary parts of $A$ and $B$ decouple from the real parts, and are neutrally stable. We are left with the following linear equation for the real parts: 
\eqh{
\left[
\begin{array}{c}
\dot{A}\\
\dot{B}
\end{array}
\right]= M(A_{\rm st},B_{\rm st}) \left[
\begin{array}{c}
{A}\\
{B}
\end{array}
\right],\quad M(A_{\rm st},B_{\rm st})=\lambda \left(
\begin{array}{cc}
1-\frac{3c+d}{c+d}&  \frac{2d}{c+d}\\
{}&{}\\
 \frac{2d}{c+d}&  1-\frac{3c+d}{c+d}
\end{array}
\right).
}
The eigenvalues of  $M(A_{\rm st},B_{\rm st})$ are equal to $\xi_1=-2$, $\xi_2=2\frac{d-c}{c+d}$, and so, the steady state is stable if $c<d$. This, together with the condition $c+d<0$ implies that the system \eqref{sysAB} possesses  a stable steady state provided $c<-|d|$ as assumed in \eqref{super2D}. Otherwise, the steady state is unstable. $\Box$


\subsection{Numerical tests for the Hookean potential}

We now compute the values of parameters $c$ and $d$ \eqref{cd} for various values of parameters $\alpha$ and $\beta$ corresponding to the slightly unstable case (close to the instability threshold). For simplicity we consider the case of unit periodic box, i.e. $L_1=L_2=L=\frac{1}{2}$, so that \eqref{fourV} gives
\eqh{
\hat{\tilde V}_{k_1,k_2}&=\frac{2\pi R^4}{z_{k_1,k_2}^2}\lr{\frac{\pi\alpha}{2}\left[J_1(z_{k_1,k_2})H_0(z_{k_1,k_2})-J_0(z_{k_1,k_2})H_1(z_{k_1,k_2})\right]-J_2(z_{k_1,k_2})},
}
where $z_{k_1,k_2}=2\pi R\sqrt{j^2+k^2}$, $\alpha=\frac{l_0}{R}$. Since we are in the periodic box, we know from Proposition \ref{Prop:Lab} that the instability appears for larger values of parameter $\beta$ than in the whole space case, i.e. for $\beta>\beta_c=\frac{24}{3-4\alpha}$.

The assumptions of Theorem \ref{Th2D} are met if
\eqh{
1+\gamma\hat{\tilde V}_{1,0}=1+\frac{\beta}{(2\pi R)^2}\lr{\frac{\pi\alpha}{2}\left[J_1(2\pi R)H_0(2\pi R)-J_0(2\pi R)H_1(2\pi R)\right]-J_2(2\pi R)}<0,
}
and
\eqh{
&1+\gamma\hat{\tilde V}_{1,1}\\
&=1+\frac{\beta}{2(2\pi R)^2}\lr{\frac{\pi\alpha}{2}\left[J_1(2\sqrt{2}\pi R)H_0(2\sqrt{2}\pi R)-J_0(2\sqrt{2}\pi R)H_1(2\sqrt{2}\pi R)\right]-J_2(2\sqrt{2}\pi R)}>0.
}
Note, that according to the definition of function $F^{\alpha,\beta}$ \eqref{defF} the above conditions are equivalent to
\[F^{\alpha,\beta}(2\pi R)<0,\quad F^{\alpha,\beta}(2\sqrt{2}\pi R)>0,\]
and from the proof of Proposition \ref{Prop:abL} we know that the rest of the eigenvalues in the assumption of Theorem \ref{Th2D} will have a good sign as well. 

We will now present computations of coefficients $c$ and $d$ defined in \eqref{cd}, that are used in Theorem \ref{Th2D} to determine the condition for the type of bifurcation \eqref{super2D}. To this purpose we choose parameter $\alpha$ in the unstable regime, here $\alpha=\frac{1}{2}$ and for several values of $R\leq L=\frac{1}{2}$ we first find the critical value of parameter $\beta$, for which the bifurcation occurs. Having this parameter we compute $c$ and $d$ using the expressions \eqref{cd} in which we take $\gamma=\frac{\beta_c}{2\pi R^4}$, we have:
\begin{enumerate}
\item for $R=\frac{1}{2}$ we have: $\beta_c=83.044$, $c=-26.327$, $d= -8.078$,
\item for $R=\frac{1}{4}$ we have: $\beta_c=31.056$, $c=-7.948$, $d=239.936$,
\item for $R=\frac{1}{8}$ we have: $\beta_c=25.544$, $c=71.726$, $d=1Â 201.065$.
\end{enumerate}
Criterion \eqref{super2D} yields that in the first case the transition is continuous while in the two following cases it is discontinuous.
Our computations are in line with the analysis in \cite{ChPa}, according to which for short range potentials (when $R/L$ is small), the transition tends to become discontinuous (first order), which corresponds to the subcritical dynamical scenario. Note that the present bifurcation analysis provides a precise criterion for the boundary
between the first order/subcritical and second order/supercritical cases.

\bigskip

\noindent {\bf Aknowledgements:} P.D. acknowledges support from the National Science Foundation (NSF) under
grants DMS-1515592 and RNMS11-07444 (KI-Net), the Engineering and Physical Sciences
Research Council (EPSRC) under grant ref. EP/M006883/1. He is on leave from CNRS, Institut de Mathematiques, Toulouse, France. He acknowledges support from the Royal Society and the Wolfson foundation through a Royal Society Wolfson Research Merit Award.
J.B. thanks the Department of Mathematics at Imperial College for hospitality, under a joint CNRS-Imperial College fellowship. E.Z. was supported by the the Department of Mathematics, Imperial College, through the Chapman Fellowship, she wishes to thank Jos\'e Antonio Carrillo for suggesting the literature and for stimulating discussions on the subject.

\footnotesize

\begin{thebibliography}{10}

\bibitem{BaSu16}
J.~W. Barrett and E.~S{{\"u}}li.
\newblock Existence of global weak solutions to compressible isentropic
  finitely extensible nonlinear bead--spring chain models for dilute polymers:
  {T}he two-dimensional case.
\newblock {\em J. Differential Equations}, 261(1):592--626, 2016.

\bibitem{BeTo11}
A.~J. Bernoff and C.~M. Topaz.
\newblock A primer of swarm equilibria.
\newblock {\em SIAM J. Appl. Dyn. Syst.}, 10(1):212--250, 2011.

\bibitem{BeCaLa09}
A.~L. Bertozzi, J.~A. Carrillo, and T.~Laurent.
\newblock Blow-up in multidimensional aggregation equations with mildly
  singular interaction kernels.
\newblock {\em Nonlinearity}, 22(3):683--710, 2009.

\bibitem{BrDeYa}
C.~P. Broedersz, M.~Depken, N.~Y. Yao, M.~R. Pollak, D.~A. Weitz, and F.~C.
  MacKintosh.
\newblock Cross-link-governed dynamics of biopolymer networks.
\newblock {\em Phys. Rev. Lett.}, 105:238101, 2010.

\bibitem{BuCl07}
G.~A. Buxton and N.~Clarke.
\newblock ``{B}ending to stretching'' transition in disordered networks.
\newblock {\em Physical review letters}, 98(23):238103, 2007.

\bibitem{CaCaPa}
J.~A. Ca{\~n}izo, J.~A. Carrillo, and F.~S. Patacchini.
\newblock Existence of compactly supported global minimisers for the
  interaction energy.
\newblock {\em Arch. Ration. Mech. Anal.}, 217(3):1197--1217, 2015.

\bibitem{CaCaSc}
J.~A. Ca{\~n}izo, J.~A. Carrillo, and M.~E. Schonbek.
\newblock Decay rates for a class of diffusive-dominated interaction equations.
\newblock {\em J. Math. Anal. Appl.}, 389(1):541--557, 2012.

\bibitem{CaChHu}
J.~A. Carrillo, A.~Chertock, and Y.~Huang.
\newblock A finite-volume method for nonlinear nonlocal equations with a
  gradient flow structure.
\newblock {\em Commun. Comput. Phys.}, 17(1):233--258, 2015.

\bibitem{CaDeMe}
J.~A. Carrillo, M.~G. Delgadino, and A.~Mellet.
\newblock Regularity of {L}ocal {M}inimizers of the {I}nteraction {E}nergy
  {V}ia {O}bstacle {P}roblems.
\newblock {\em Comm. Math. Phys.}, 343(3):747--781, 2016.

\bibitem{CaMcVi03}
J.~A. Carrillo, R.~J. McCann, and C.~Villani.
\newblock Kinetic equilibration rates for granular media and related equations:
  entropy dissipation and mass transportation estimates.
\newblock {\em Rev. Mat. Iberoamericana}, 19(3):971--1018, 2003.

\bibitem{ChPa}
L.~Chayes and V.~Panferov.
\newblock The {M}c{K}ean-{V}lasov equation in finite volume.
\newblock {\em J. Stat. Phys.}, 138(1-3):351--380, 2010.

\bibitem{DeDePe}
P.~Degond, F.~Delebecque, and D.~Peurichard.
\newblock Continuum model for linked fibers with alignment interactions.
\newblock {\em Math. Models Methods Appl. Sci.}, 26(2):269--318, 2016.

\bibitem{DeLiRi}
P.~Degond, J.-G. Liu, and C.~Ringhofer.
\newblock Evolution of wealth in a non-conservative economy driven by local
  {N}ash equilibria.
\newblock {\em Philos. Trans. R. Soc. Lond. Ser. A Math. Phys. Eng. Sci.},
  372(2028):20130394, 15, 2014.

\bibitem{BeCh06}
M.R. D'Orsogna, Y.L. Chuang, A.L. Bertozzi, and L.S. Chayes.
\newblock Self-propelled particles with soft-core interactions: patterns,
  stability, and collapse.
\newblock {\em Phys Rev Lett.}, 96(10)(104302), Mar 17 2006.

\bibitem{FiRu66}
M.~E. Fisher and D.~Ruelle.
\newblock The stability of many-particle systems.
\newblock {\em J. Mathematical Phys.}, 7:260--270, 1966.

\bibitem{Gardel2004}
M.~L. Gardel, J.~H. Shin, F.~C. MacKintosh, L.~Mahadevan, P.~Matsudaira, and
  D.~A. Weitz.
\newblock Elastic behavior of cross-linked and bundled actin networks.
\newblock {\em Science}, 304(5675):1301--1305, 2004.

\bibitem{HaIo11}
M.~Haragus and G.~Iooss.
\newblock {\em Local bifurcations, center manifolds, and normal forms in
  infinite-dimensional dynamical systems}.
\newblock Universitext. Springer-Verlag London, Ltd., London; EDP Sciences, Les
  Ulis, 2011.

\bibitem{KolCarBer}
T.~Kolokolnikov, J.~A. Carrillo, A.~Bertozzi, R.~Fetecau, and M.~Lewis.
\newblock Emergent behaviour in multi-particle systems with non-local
  interactions [{E}ditorial].
\newblock {\em Phys. D}, 260:1--4, 2013.

\bibitem{McK67}
H.~P. McKean, Jr.
\newblock Propagation of chaos for a class of non-linear parabolic equations.
\newblock In {\em Stochastic {D}ifferential {E}quations ({L}ecture {S}eries in
  {D}ifferential {E}quations, {S}ession 7, {C}atholic {U}niv., 1967)}, pages
  41--57. Air Force Office Sci. Res., Arlington, Va., 1967.

\bibitem{MoEd99}
A.~Mogilner and L.~Edelstein-Keshet.
\newblock A non-local model for a swarm.
\newblock {\em J. Math. Biol.}, 38(6):534--570, 1999.

\bibitem{PeDeDe}
D.~Peurichard, F.~Delebecque, A.~Lorsignol, C.~Barreau, J.~Rouquette,
  X.~Descombes, L.~Casteilla, and P.~Degond.
\newblock Simple mechanical cues could explain adipose tissue morphology.
\newblock {\em submitted}.

\bibitem{Ru69}
D.~Ruelle.
\newblock {\em Statistical mechanics: {R}igorous results}.
\newblock W. A. Benjamin, Inc., New York-Amsterdam, 1969.

\bibitem{SiSlTo15}
R.~Simione, D.~Slep{\v{c}}ev, and I.~Topaloglu.
\newblock Existence of ground states of nonlocal-interaction energies.
\newblock {\em J. Stat. Phys.}, 159(4):972--986, 2015.

\end{thebibliography}

\end{document}